\documentclass[11pt]{article}

\usepackage{amsmath,amsmath,amsthm,amscd,amssymb,verbatim,epsf}

\numberwithin{equation}{section}

\usepackage{amsthm, amssymb, amscd}

\newtheorem{pro}{Proposition}[section]
\newtheorem{thm}[pro]{Theorem}
\newtheorem{lem}[pro]{Lemma}

\newtheorem{cor}[pro]{Corollary}

\newtheorem{defi}[pro]{Definition}

\newtheorem{adj}[pro]{Adjustment}
\newtheorem{note}[pro]{Note}
\newtheorem{cond}[pro]{Condition}
\newtheorem{nota}[pro]{Notation}

\def\G{{\Gamma}}
 \def\d{{\delta}}
 
 \def\e{{\epsilon}}
 
 \def\L{{\Lambda}}

  \def\O{{\Omega}}
   \def\s{{\sigma}}
 
 \def\a{{\alpha}}
 \def\b{{\beta}}
 \def\p{{\partial}}
 
 \def\ra{{\rightarrow}}
 \def\lra{{\longrightarrow}}
 \def\g{{\gamma}}

 \def\c{{\mathbb C}}
 
 \def\z{{\mathbb Z}}
 \def\2{{\mathbb Z_2}}
 
 \def\t{{\tau}}

 \def\sl2{{SL(2,\mathbb C)}}
 
 \def\qed{{\hspace{2mm}{\small $\diamondsuit$}}}

 \def\pf{{\noindent{\bf Proof.\hspace{2mm}}}}

 \def\sl{{{\mbox{\tiny $\L$}}}}

\def\lra{\longrightarrow}

\def\d{\delta}

\def\a{\alpha}
\def\b{\beta}
\def\t{\tau}
\def\L{\Lambda}
\def\e{\epsilon}

\def\g{\gamma}

\def\s{\sigma}
\begin{document}

\begin{center}
{\bf  Quasi-Fuchsian Surfaces In  Hyperbolic Link Complements}
\end{center}

\begin{center}
 Joseph D. Masters and Xingru Zhang
\end{center}

{\bf Abstract.} We show that every hyperbolic link complement contains closed
quasi-Fuchsian surfaces. As a consequence, we obtain the result that on a
hyperbolic link complement, if we remove  from each cusp of the manifold a
certain finite set of  slopes, then all remaining Dehn fillings on the link
complement yield  manifolds with closed immersed incompressible surfaces.

\vspace{10mm} \section{Introduction}

 By a {\it link complement} we mean, in this
paper, the complement of a link in a closed connected orientable $3$-manifold. A
link complement is said to be hyperbolic if it admits a complete hyperbolic
metric of finite volume.
 By a \textit{surface} we mean, in this paper,
 the complement of a finite (possibly empty) set of points in
 the interior of  a compact, orientable $2$-manifold (which may
 not be connected).
By a \textit{surface in a $3$-manifold} $W$, we mean a continuous, proper map $f:S \ra W$ from a surface
$S$ into $W$.
  A surface $f:S\ra W$ in a 3-manifold $W$ is said to be
 \textit{connected} if and only if $S$ is connected.
  A surface $f:S\ra W$ in a 3-manifold $W$ is said to be
 \textit{incompressible}
 if   each component $S_j$ of $S$ is not a $2$-sphere and the induced homomorphism
 $f^*: \pi_1(S_j, s)\ra \pi_1(W, f(s))$ is injective for  any
 choice of base  point $s$ in $S_j$.
 A surface $f:S\ra W$ in a 3-manifold $W$
 is said to be {\it essential} if it is incompressible and
 for each component $S_j$ of $S$, the map $f: S_j\ra W$ cannot be properly
  homotoped  into a boundary component or an end component
 of $W$.

 Connected essential surfaces in  hyperbolic link complements can be divided
into three mutually exclusive geometric types: quasi-Fuchsian surfaces,
geometrically infinite surfaces, and essential surfaces with accidental
parabolics. Geometrically these three types of surfaces can be characterized by
their limit sets as follows: a connected essential surface $f:S\ra M$ in a
hyperbolic link complement $M$ is Quasi-Fuchsian if and only if the limit set of
the subgroup $f^*(\pi_1(S))\subset \pi_1(M)$ is a Jordon circle in the boundary
$2$-sphere of the hyperbolic $3$-space $\mathbb H^3$; is geometrically infinite
if and only if the limit set of $f^*(\pi_1(S))$ is the whole $2$-sphere; and is
having  accidental parabolics otherwise. Topologically these three types of
surfaces can be characterized
 as follows: a connected essential surface $f:S\ra M$ in a
hyperbolic link complement $M$   is geometrically infinite if and only if
  it can be lifted (up to homotopy) to a fiber in some finite
 cover of $M$; is having accidental parabolics
  if and only  $S$ contains a closed curve which cannot be freely homotoped
  in $S$ into a cusp of $S$ but can be freely homotoped in $M$
  into a cusp of $M$,
  and is quasi-Fuchsian otherwise.

In \cite{MZ} it was shown that every hyperbolic knot complement
contains closed quasi-Fuchsian surfaces.
In this paper we extend this  result to hyperbolic link complements.

\begin{thm}\label{main}
Every hyperbolic link complement contains  closed
quasi-Fuchsian surfaces. \end{thm}

This yields directly the following  consequence.

\begin{cor}\label{cor}
For every given hyperbolic link complement $M$, if we remove certain
 finitely many  slopes from each
cusp of $M$, then all  remaining Dehn fillings  produce manifolds which  contain closed
incompressible
surfaces.
\end{cor}

 We note that Corollary \ref{cor} would also be a consequence of Khan and Markovic's recent
 claim that every closed hyperbolic 3-manifold contains a surface subgroup.

 This paper is an extension of \cite{MZ} where the existence of closed quasi-Fuchsian
surfaces in any hyperbolic knot complement was proved. The proof of Theorem \ref{main} follows
essentially the approach given in \cite{MZ}. To avoid repetition, we shall assume the reader is
familiar with
 the machinery laid out in \cite{MZ}.
 In particular we shall use most of the terms and properties about hyperbolic $3$-manifolds
 and about groups  established in \cite{MZ}, without recalling  them in detail, and shall omit
 details of constructions and proof of assertions whenever they are
 natural generalization of counterparts of \cite{MZ}.

 To help the reader to get a general idea about which parts of our
 early arguments are needed to be adjusted  nontrivially, we first very briefly
 recall how a closed quasi-Fuchsian surface was  constructed in a hyperbolic
knot complement.
  We started with a pair of connected bounded embedded
quasi-Fuchsian surfaces in a given hyperbolic knot exterior $M^-$
(which is truncation of a hyperbolic knot complement $M$) with distinct
boundary slopes in $\p  M^-$.
We then considered  two hyperbolic  convex $I$-bundles resulting from the two
corresponding quasi-Fuchsian surface groups. By a careful ``convex  gluing'' of
two suitable finite covers of  some truncated versions of the  two
 $I$-bundles, and then ``capping off convexly'' by a solid cusp, we
constructed a convex hyperbolic $3$-manifold $Y$ with a local isometry $f$ into
the given hyperbolic knot complement  $M$. The manifold $Y$ had non-empty boundary
each component of which provided a closed quasi-Fuchsian surface in $M$ under
the map $f$. To find the required finite covers of the truncated $I$-bundles and
at the same time to lift certain immersions  to embeddings, we needed a stronger
version of subgroup separability property for surface groups with boundary, which was
proved using Stallings' folding
graph techniques.

Now to extend the construction to work for  hyperbolic link complements, we first  need to prove,  for
any given hyperbolic link exterior $M^-$, the existence of two properly embedded bounded quasi-Fuchsian
surfaces $S_i^-$, $i=1,2$, in $M^-$, each of which is not necessarily connected, with the crucial
property that for each of $i=1,2$ and each component $T_j$ of $\p M^-$, $S_i^-\cap T_j$ is a non-empty
set of simple closed essential curves, and furthermore the slope of the curves $S_1^-\cap T_j$ is
different from that of $S_2^-\cap T_j$ for each $T_j$. The proof
 of this result,  given in Section \ref{bs}, is based on work of
Culler-Shalen \cite{CS} and  Cooper-Long \cite{CL} and Thurston \cite{T},
 making use of the $SL_2(\c)$ character
variety of the link exterior $M^-$ and some special properties  of essential surfaces in hyperbolic
$3$-manifolds with accidental parabolics. With the given two surfaces $S_i^-$, $i=1,2$, we may construct
two corresponding convex $I$-bundles (in the current situation each $I$-bundle may not be connected).
Following the approach in \cite{MZ} we still
 want to choose a suitable cover for each component of each  of the truncated $I$-bundles, and ``convexly
glue'' them all together in certain way and ``convexly cap off'' with $m$ (which is the number of
components of $\p  M^-$) solid cusps, to form a convex hyperbolic $3$-manifold $Y$, with a local isometry
$f$ into $M$, such that the boundary of $Y$ is a non-empty set, each component of which is mapped
 by $f$ to a closed quasi-Fuchsian surface in $M$.

 As before, we want to choose the cover so that the boundary components of $S_i^-$
 unwrap as much as possible.  And if the two surfaces $S_i^-$ are connected, our previous arguments
 go through with very little change.
 However, if the surfaces $S_i^-$ are disconnected, complications arise.  In order to piece the
 different covers together, we need to know that they all have the same degree.  And this turns out to
 require a non-trivial strengthening of our previous separability result; see Theorem \ref{each large n}.
  The proof of this  property uses a careful refinement of the folding graph arguments used in \cite{MZ}.

\section{Cusped qausi-Fuchsian surfaces in hyperbolic link
complements}\label{bs}

From now on let $M$ be a given hyperbolic link complement of $m\geq 2$ cusps. For each of $i=1,...,m$,
let $C_i$ be a fixed $i$-th cusp of $M$ which is geometric, embedded and small enough so that
$C_1,...,C_m$ are mutually disjoint. The complement of the interior of $C_1\cup...\cup C_m$
 in $M$,
which we denote by $M^-$, is a compact, connected and orientable $3$-manifold whose boundary is a set of
$m$ tori. We call $M^-$ a \textit{truncation} of $M$. Let $T_k=\p C_k$, $k=1,...,m$.
 Then $\p
M^-=T_1\cup...\cup T_m$.

\begin{lem}\label{slopes}
There are two embedded essential
 quasi-Fuchsian surfaces $S_1$ and $S_2$ in $M$  (each $S_i$ may not be
 connected) such that for  each of $i=1,2$ and  each of $k=1,...,m$,
 $S_i\cap T_k$
 is a nonempty set of parallel simple closed essential
curves in $T_k$ of slope $\lambda_{i,k}$ and $\lambda_{1,k}\ne \lambda_{2,k}$.
\end{lem}

\pf It is equivalent to show that
 the truncation $M^-$ of $M$
contains two properly  embedded bounded essential
  surfaces $S_1^-$ and $S_2^-$
  such that:
  \newline
   (i) For each of $i=1,2$, each component of $S_i^-$ is not a fiber or semi-fiber
   of $M^-$.
 \newline
   (ii) For each of $i=1,2$, any closed curve in $S_i^-$
  that can be freely homotoped in $M^-$ into $\p M^-$
 can also be freely homotoped in $S_i^-$ into $\p S_i^-$.
 \newline
  (iii) For each of  $i=1,2$  and each of $k=1,...,m$,
  $S_i^-$ has non-empty boundary on $T_k$ of boundary slope
   $\lambda_{i,k}$ and $\lambda_{1,k}\ne \lambda_{2,k}$.

Let $\{\g_k\subset T_k; k=1,...,m\}$ be any given set of $n$ slopes.
By \cite[Theorem 3]{CS}, there is
a properly embedded essential surface $S_1^-$ (maybe disconnected) in $M^-$ with the following
properties (in fact the surface $S_1^-$ is obtained through a nontrivial group action on a simplicial
tree associated to an ideal point of a curve in a component of the $SL(2,\c)$-character variety of $M^-$
which contains the character of a discrete faithful representation of $\pi_1(M^-)$):

\noindent (1) No component of $S_1^-$ is a fiber or semi-fiber of $M^-$.

\noindent (2) For each of $k=1,...,m$, $S_1^-$ has non-empty boundary on $T_k$ of boundary slope
$\lambda_{1,k}$ which is different from $\g_k$.

\noindent (3) If  an element of $\pi_1(M^-)$ is freely homotopic to a
curve in $M^-\setminus S_1^-$, then it is contained in a vertex stabilizer of
the action on the tree.

\noindent (4) If  an element of $\pi_1(M^-)$ is freely homotopic to
$\g_k$, then it is not contained in any vertex stabilizer of the action on the
tree and thus must intersect $S_1^-$.

\noindent (5) If  an element of $\pi_1(M^-)$ is freely homotopic to a
curve in $S_1^-$, then it  is  contained in an edge  stabilizer of the tree.
\newline
It follows that

\noindent (6) If an element of $\pi_1(M^-)$ is freely homotopic to a
simple closed essential curve in $T_k$ whose slope is different from
$\lambda_{1,k}$, then it is not contained in any vertex stabilizer of the action
on the tree.

 Let $S^-_{1,j}$, $j=1,...,n_1$, be the components of $S^-_1$. If some
$S^-_{1,j}$ has a
   closed curve  which  cannot be freely homotoped in $S_1^-$ into $\p S_1^-$
  but can be freely homotoped in $M^-$ into $\p M$,
  then arguing as in   \cite[Lemma 2.1]{CL}, we see that there is an
   embedded annulus $A$ in $M^-\setminus S^-_{1,j}$ such that one
   boundary component, denoted $a_1$,  of $A$
lies in $S_{1,j}^-$ and  is not  boundary parallel in $S^-_{1,j}$, and the other
boundary component, denoted $a_2$, of $A$ is contained in some boundary
component $T_k$ of $M^-$. By Properties  (5) and (6) listed above, we have

\noindent  (7) $a_2\subset T_k$ must have the slope $\lambda_{1,k}$.

Now consider in  $A$ the intersection   set $A\cap (S^-_1-S^-_{1,j})$ of $A$
with other components of $S_1^-$. By Property (7), we may assume that $\p A\cap
(S^-_1-S^-_{1,j})=\emptyset$. Thus by  proper isotopy of $(S_1^--S^-_{1,j}, \p
(S^-_1-S^-_{1,j}))\subset (M^-,\p M^-)$ and surgery (if necessary) we may assume
that each component of $A\cap (S^-_1-S^-_{1,j})$ is a circle which is isotopic
in $A$ to the center circle of $A$ and if the component is contained in
$S^-_{1,j'}$, then it is not boundary parallel in $S^-_{1,j'}$. So the component
of $A\cap (S^-_1-S^-_{1,j})$, denoted $a_1'$,  which is closest to $a_2$ in $A$,
 cuts out from $A$ an sub-annulus $A'$ which is properly embedded in $M^-\setminus S^-_1$ such that
 $a_1'$ lies in $S^-_{1,j'}$, for some $j'$, and is not boundary
 parallel in $S^-_{1,j'}$. So we may perform the annulus compression on $S^-_{1,j'}$
 along $A'$ to get an essential surface
 which still satisfies the properties (1)-(6) above (because the new resulting
  surface can be considered as
 a subsurface of the old surface $S_1^-$ and because of property (7)) but has
 larger Euler characteristic.
 Thus such annulus compression must terminate in a finite number of times.
So eventually we end up with a surface,  which we still denote by $S_1^-$,  satisfying the condition

\noindent (8) Any closed curve in $S_1^-$
  that can be freely homotoped in $M^-$ into $\p M$
 can  be freely homotoped in $S_1^-$ into $\p S_1^-$.

Now letting $\g_k=\lambda_{1,k}, k=1,...,m$,
and repeating the above arguments, we may get another
properly embedded essential surface $S_2^-$ such that

\noindent (1') Each component of $S_2^-$ is not a fiber or semi-fiber of $M^-$.

\noindent (2') For each of $k=1,...,m$, $S_2^-$ has non-empty boundary
on $T_k$ of boundary slope
$\lambda_{2,k}$ which is different from $\lambda_{1,k}$.

\noindent (8') Any closed curve in $S_2^-$
  that can be freely homotoped in $M^-$ into $\p M$
 can be freely homotoped in $S_2^-$ into $\p S_2^-$.

\noindent So $S_1^-$ and $S_2^-$ satisfy conditions (i), (ii) and (iii) listed
above. The lemma is thus proved. \qed

Let $S_i,i=1,2$  be the two surfaces provided by Lemma \ref{slopes}.
By taking disjoint parallel copies
of some components of $S_i$ (if necessary), we may and shall assume

\begin{cond}\label{at least two}
{\rm For each $i=1,2$ and $k=1,...,m$, $S_i\cap T_k$ has a positive, even number of  components}.
\end{cond}

\begin{nota}\label{dik}
{\rm  Let $S_{i,j}$, $j=1,...,n_i$, be components of $S_i$, $i=1,2$. Let $i_*$ be the number such that
$\{i,i_*\}=\{1,2\}$ for $i=1,2$.
 Let $S_i^-=S_i\cap M^-$ and let $\p_k
S_{i,j}^-$ be the boundary components of $S^-_{i,j}$ on $T_k$ (which may be empty for some $j$'s) and
let $\p_k S_i^-=\cup_j \p_k S_{i,j}^-$. Now for each $i=1,2, k=1,...,m$, let $d_{i,k}$ be the geometric
intersection number in $T_k$ between a  component of $\p_k S_i^-$  and the whole set
 $\p_k S_{i_*}$.
 Obviously  $d_{i, k}$ is independent of the choice of the component of $\p_k S_i^-$.
By Condition \ref{at least two}, $d_{i,k}\geq 2$ is even for each $i,k$. Now set
 $$ d_i=lcm\{d_{i,k}; \;k=1,...,m\},$$ the
(positive) least common multiple. Then $d_i\geq 2$ is even for each $i=1,2$.}
\end{nota}

 Let $\mathbb H^3$ be the
hyperbolic $3$-space in the upper half space model, let $S_{\infty}^2$ be the $2$-sphere at $\infty$  of
$\mathbb H^3$ and let $\overline{\mathbb H}^3 =\mathbb H^3\cup S_{\infty}^2$.

By Mostow-Prasad  rigidity,  the fundamental group of $M$ (for any fixed choice of base point) can be
uniquely identified as a discrete torsion free subgroup $\G$ of $Isom^+(\mathbb H^3)$ up to conjugation
in $Isom(\mathbb H^3)$ so that $M=\mathbb H^3/\G$. We shall fix one such identification. Let $p:\mathbb
H^3\ra M$ be the corresponding covering map.

For the given  surface $S_{i,j}$ in $M$ (for each $i, j$), we identify its fundamental group with a quasi-Fuchsian subgroup
$\G_{i,j}$ of $\G$ as follows. As $S_{i, j}$ is embedded in $M$ we may consider it as a submanifold of
$M$. Fix a component $\tilde S_{i,j}$ of $p^{-1}(S_{i,j})$ (topologically $\tilde S_{i,j}$ is an  open
disk in $\mathbb H^3$), there is a subgroup $\G_{i,j}$ in the stabilizer
 of $\tilde S_{i,j}$ in $\G$ such that
 $S_{i,j}=\tilde S_{i,j}/\G_{i,j}$.

 Note that the limit set $\L_{i,j}$ of $\G_{i,j}$ is
 a Jordan circle in the $2$-sphere $S^2_{\infty}$
 at the $\infty$ of $\mathbb H^3$.
 Let $H_{i,j}$ be the convex hull of
 $\L_{i,j}$ in $\mathbb H^3$.

 Let ${\cal B}_k=p^{-1}(C_k)$, $k=1,...,m$, and ${\cal B}=p^{-1}(C)$.
 Then by our assumption on $C$,  ${\cal B}$ is a set of mutually disjoint horoballs
 in $\mathbb H^3$.
 Let $B$ be a component of ${\cal B}$ and let $\p  B$ be the frontier of  $B$ in $\mathbb H^3$.
 Then $\p B$ with the
 induced metric is isometric to a  Euclidean plane.
 We shall simply call $\p B$ a Euclidean plane.
  A strip between two parallel Euclidean lines in $\p B$
 will be called a Euclidean strip in $\p B$.
 Note that every Euclidean line in $\p B$ bounds a totally geodesic
 half plane in $B$ (which is perpendicular to $\p B$).
  By a {\it $3$-dimensional strip region} in $B$ we mean
 a region in $B$ between two
 totally geodesic half planes in $B$ bounded by two parallel disjoint Euclidean lines
 in $\p B$.

 \begin{lem}\label{strip}
 If the cusp set $C=C_1\cup ...\cup C_m$ of $M$ is small enough, then
 for each component  $B$ of ${\cal B}$ whose point at $\infty$  is a parabolic fixed point
 of $\G_{i,j}$,
  $H_{i,j}\cap B$ is a $3$-dimensional strip region in $B$.
  \end{lem}

 \pf The proof is similar to that of \cite[Lemma 5.2]{MZ}.
 \qed

From now on we assume that $C$ has been chosen so that Lemma \ref{strip} holds for all $i=1,2,
j=1,...,n_i$. For a fixed small $\e>0$, let $X_{i,j}$ be the $\e$-collared neighborhood of $H_{i,j}$ in
$\mathbb H^3$. Then it follows from Lemma \ref{strip} that for each component  $B$ of ${\cal B}$ whose
point at $\infty$  is a parabolic fixed point
 of $\G_{i,j}$,,
$X_{i,j}\cap B$ is a $3$-dimensional strip region in $B$,
 for all $i=1,2, j=1,..., n_i$, by geometrically shrinking $C$ further if necessary.

Note that $X_{i,j}$  is a   metrically complete and strictly convex
hyperbolic $3$-submanifold of $\mathbb{H}^3$ with $C^1$ boundary,
 invariant under the action of $\G_{i,j}$. Let $${\cal
B}_{i,j}=\{X_{i,j}\cap B; \mbox{$B$ a component of $ {\cal B}$ based at a
parabolic fixed point of $\G_{i,j}$}\}.$$ We call ${\cal B}_{i,j}$ the
\textit{horoball region} of $X_{i,j}$. Let $X_{i,j}^-=X_{i,j}\setminus {\cal B}_{i,j}$, and
call  $X_{i,j}^-\cap \p {\cal B}_{i,j}$ the \textit{parabolic boundary} of $X_{i,j}^-$,
denoted by $\p_p X_{i,j}^-$. Note that $ X_{i,j}^-$ is locally
convex everywhere except on its parabolic boundary.

Each of $X_{i,j}$,
 ${\cal B}_{i,j}$, $X_{i,j}^-$ and $\p_p X_{i,j}^-$ is invariant under the
action of $\G_{i,j}$. Let $Y_{i,j}=X_{i,j}/\G_{i,j}$, which is a metrically
complete and strictly convex hyperbolic $3$-manifold with
boundary. Topologically $Y_{i,j}=S_{i,j}\times I$,
 where $I = [-1, 1]$. There is a
local isometry $f_{i,j}$ of $Y_{i,j}$ into $M$, which is induced from the
covering map $\mathbb H^3/\G_{i,j}\;\lra \;M$ by restriction on $Y_{i,j}$,
 since $Y_{i,j}=X_{i,j}/\G_{i,j}$ is a submanifold of $\mathbb H^3/\G_{i,j}$. Also
$p|_{X_{i,j}}=f_{i,j}\circ p_{i,j}$, where $p_{i,j}$ is the universal covering map
$X_{i,j}\ra Y_{i,j}=X_{i,j}/\G_{i,j}$.
 Let $Y_{i,j}^-=X_{i,j}^-/\G_{i,j}$,  let ${\cal C}_{i,j}={\cal B}_{i,j}/\G_{i,j}$,
and let  $\p_p Y_{i,j}^-=\p_p X_{i,j}^-/\G_{i,j}$. We call ${\cal C}_{i,j}$ the cusp part of $Y_{i,j}$,
and call $\p_pY_{i,j}^-$ the parabolic boundary of $Y_{i,j}^-$, which is the frontier of $Y_{i,j}^-$ in
$Y_{i,j}$ and is also the frontier of ${\cal C}_{i,j}$ in $Y_{i,j}$. Each component of $\p_p Y_{i,j}^-$
is a Euclidean annulus. The manifold $Y_{i,j}^-$ is locally convex everywhere except on its parabolic
boundary. Topologically $Y_{i,j}^-=S_{i,j}^-\times I$.

As in \cite[Section 5]{MZ}, we fix a product $I$-bundle structure for $Y_{i,j}=S_{i,j}\times I$ such
that each component of ${\cal C}_{i,j}$ has the induced $I$-bundle structure which is the product
 of a totally geodesic cusp annulus and the $I$-fiber (i.e. we assume that
 $(S_{i,j}\times \{0\})\cap {\cal C}_{i,j}$
 is a set of totally geodesic cusp annuli).
We let every free cover of $Y_{i,j}$ have the induced $I$-bundle structure.
 In particular $X_{i,j}$ has the induced $I$-bundle structure from that of $Y_{i,j}$,
 and this structure is preserved by the action of $\G_{i,j}$; i.e. every
 element of $\G_{i,j}$ sends an $I$-fiber of $X_{i,j}$ to an $I$-fiber of $X_{i,j}$.
Similar to  \cite[Corollary 5.6]{MZ}, we have

\begin{lem}\label{length}
For each of $i=1,2, j=1,...,n_i$, there is an upper bound for the lengths of the
$I$-fibers of $X_{i,j}$. \qed
\end{lem}

 The restriction of the map $f_{i,j}$ on the center surface
 $S_{i,j}\times \{0\}$ of $Y_{i,j}=S_{i,j}\times I$
 may not be an embedding in general
 but by  Lemma \ref{strip} we may and shall assume that the map is an embedding
when restricted on  $(S_{i,j}\times \{0\})\cap {\cal
C}_{i,j}$.
 We now replace our original embedded surface $S_{i,j}$
 by the center surface $f_{i,j}:S_{i,j}\times\{0\}$ and
 we simply denote $S_{i,j}\times \{0\}$ by $S_{i,j}$.

The restriction map  $f_{i,j}:(Y_{i,j}^-,  \p_p Y_{i,j}^-)\ra (M^-,\p M^-)$ is
a proper map of pairs and  $f_{i,j}:
 (S_{i,j}^-,\p S_{i,j}^-)\ra (M^-, \p M^-)$ is a
 proper map which is an embedding on $\p S_{i,j}^-$ (This property will remain
 valid if we shrink the cusp $C$ of $M$ geometrically).
 In fact $f_{i,j}(\p S_{i,j}^-)$ are embedded Euclidean circles in
 $\p M^-$.
 Hence boundary slopes of the new quasi-Fuchsian surfaces
 $f_{i,j}:
 (S_{i,j}^-,\p S_{i,j}^-)\ra (M^-, \p M^-)$
 are defined and are the same as those of the original embedded surfaces
 $S_{i,j}^-$.

\begin{note}\label{still defined}
{\rm As $f_{i, j}: \p S_{i, j}^-\ra \p M^-$ is an embedding, we sometimes
 simply consider $\p S_{i, j}^-$ as subset of $\p M^-$, for each $i=1,2,
 j=1,..., n_i$.
 By choosing  a slightly different center surface for $Y_{i, j}$ (if necessary),
 we may assume that the components of
  $\{\p S_{i, j}^-, j=1,..., n_i\}$ are mutually disjoint in $\p M^-$,
  for each fixed $i=1,2$.
  So the numbers $d_{i, k}$, $d_i$ defined in Notation \ref{dik}
  remain well defined
  for the current surface $f_{i, j}:(S_{i, j}^-, \p S_{i, j}^-)\ra (M, \p M)$
  and are the same numbers as given there, for all $i, j$.
 Also $\p_k S_{i, j}^-$ remain defined as before for all $i, j$. }\end{note}

Let $\tilde S_{i, j}$ and $\tilde S_{i,j}^-$ be the corresponding
 center surfaces of $X_{i,j}$
and $X_{i, j}^-$ respectively.

Note also that if a component of $\p S_{i,j}^-$
 intersects a component of $\p S_{i_*,j'}^-$
 in some component $T_k$ of $\p M^-$, then they intersect
  geometrically in $T_k$, and
  their intersection points in $T_k$ are one-to-one corresponding
to the geodesic rays of $f_{i,j}(S_{i,j})\cap f_{i_*,j'}(S_{i_*,j'})\cap C_k$.

 We fix an orientation for $S_{i,j}$, and let $S_{i,j}^-$ and $\p S_{i,j}^-$
 have the induced orientation.

\section{Construction of intersection pieces $K_{i,j}$}

Suppose that $\p S_{i,j}^-$
 intersects  $\p S_{i_*,j'}^-$ for some $j,j'$.
 We  construct the ``intersection pieces'' $K_{i, j,j'}$ and $K_{i_*,j',j}$
 between
 $Y_{i,j}$
and $Y_{i_*,j'}$ in a similar fashion as in \cite[Section 6]{MZ} such that
\newline
(1) $K_{i, j,j'}$ and $K_{i_*,j',j}$ are isometric.
\newline
(2) Each component of $K_{i, j,j'}$ or of $K_{i_*,j',j}$ is a metrically
 complete convex hyperbolic
$3$-manifold.
\newline
(3) There are local isometries $g_{i,j,j'}: K_{i,j,j'}\ra Y_{i,j}$ and
$g_{i_*,j',j}:K_{i_*,j',j}\ra Y_{i_*,j'}$.
\newline
(4) $K^-_{i,j,j'}$ and $K^-_{i_*,j',j}$ (which are the truncated
versions of $K_{i, j,j'}$ and $K_{i_*,j',j}$ respectively) are compact.
\newline
(5) Each component of the parabolic boundary
 $\p_p K_{i,j,j'}^-$ of $K_{i,j,j'}^-$ is a Euclidean
parallelogram, the number of cusp ends of $K_{i,j,j'}$ is precisely the
 number of intersection points
between $f_{i,j}(\p S_{i,j}^-)$ and $f_{i_*,j'}(\p S_{i_*,j'}^-)$.
Similar properties  hold for
$K_{i_*,j',j}$.
\newline
(6) The restriction of
$g_{i,j,j'}$ to $K_{i,j,j'}\setminus K_{i,j,j'}^-$ is an embedding
and so is the restriction of
$g_{i_*,j',j}$ to $K_{i_*,j',j}\setminus K_{i_*,j',j}^-$.
\newline
(7) $f_{i,j}(g_{i,j,j'}(K_{i,j,j'}\setminus K_{i,j,j'}^-))$
contains $f_{i,j}(S_{i,j})\cap f_{i_*,j'}(S_{i_*,j'})
\cap C$ (the latter  is a set of geodesic rays) and so does
$f_{i_*,j'}(g_{i_*,j',j}(K_{i_*,j',j}\setminus K_{i_*,j',j}^-))$.

Let $K_{i,j}$ be the disjoint union of these $K_{i,j,j'}$ over such $j'$. Then the number of components
of $\p_p K_{i,j}^-$ is precisely the number of intersection points between $\p S_{i,j}^-$ and $\p
S_{i_*}^-$. In fact there is a canonical one-to-one correspondence between components of $\p_p K_{i,j}^-$
and the intersection points between $\p S_{i,j}^-$ and $\p S_{i_*}^-$.

 Let $K_{i}$ be
the disjoint union of these $K_{i,j}$. Then the number of cusp ends of $K_{i}$ is precisely the number
of intersection points between $\p S_{i}^-$ and $\p S_{i_*}$ and there is an isometry between $K_1$ and
$K_2$.

\section{Construction of $J_{i,j}$, $J_{i,j}^-$ , $\hat J_{i,j}$
and $C_n(J_{i,j}^-)$}

As in \cite[Section 6]{MZ},
we fix a number  $R>0$  bigger than  the number $R(\e)$ provided
in \cite[Proposition 4.5]{MZ}  and also bigger than the upper bound
 provided by Lemma \ref{length} for
 the lengths of $I$-fibers of $X_{i,j}$ (for each of $i=1,2$, $j=1,...,n_i$).
As in \cite[Section 6]{MZ}, we define and construct
 the  {\it
abstract $R$-collared neighborhood of $K_{i,j}$ with respect to $X_{i,j}$}
which is denoted by $AN_{(R,X_{i,j})}(K_{i,j})$.
Also define the truncated version
$(AN_{(R,X_{i,j})}(K_{i,j}))^-$, the parabolic boundary $\p_p (AN_{(R,X_{i,j})}(K_{i,j}))^-$
and the cups region $AN_{(R,X_{i,j})}(K_{i,j})\setminus(AN_{(R,X_{i,j})}(K_{i,j}))^-$
accordingly.

 Now as in \cite[Section 7]{MZ}, we construct a connected metrically complete,
 convex, hyperbolic $3$-manifold
 $J_{i,j}$ with a local isometry $g_{i,j}:J_{i,j}\ra Y_{i,j}$
 such that $J_{i,j}$ contains
$AN_{(R,X_{i,j})}(K_{i,j})$ as a hyperbolic submanifold, and
 $J_{i,j}\setminus AN_{(R,X_{i,j})}(K_{i,j})$  is a compact $3$-manifold $W_{i,j}$
 (which may not be
connected). Also  $W_{i,j}$ is disjoint from $AN_{(R,X_{i,j})}(K_{i,j})\setminus
(AN_{(R,X_{i,j})}(K_{i,j}))^-$,  the parabolic boundary  $\p_p J_{i,j}^-$ of
$J_{i,j}^-$ is equal to the  parabolic boundary of
$(AN_{(R,X_{i,j})}(K_{i,j}))^-$, and $g_{i,j}|: (J_{i,j}^-, \p_p J_{i,j}^-)\ra
(Y_{i,j}^-, \p_pY_{i,j}^-)$ is a proper map of pairs.

Each component of $\p_p J_{i,j}^-$ is  a Euclidean parallelogram
and thus can be capped off  by a convex $3$-ball.
Let $\hat J_{i,j}$ be the resulting manifold after capping off
all components of $\p_p J_{i,j}^-$. Then $\hat J_{i,j}$
 is a connected, compact, convex $3$-manifold with a local
isometry (which we still denote by $g_{i,j}$) into $Y_{i,j}$.

The  number of  components of $\p_p J_{i,j}^-$ is equal to the number of
components of $\p_p K_{i,j}^-$, and the former is an abstract $R$-collared
neighborhood of the latter with respect to $\p_p X_{i,j}^-$.

\begin{note}{\rm
The components of $\p_p J_{i,j}^-$ are canonically
one-to-one correspond to
the intersection points of $\p S_{i, j}^-$ with
$\p S_{i_*}^-$.}
\end{note}

Now as  in \cite[Section 8]{MZ}, we construct, for each sufficiently large integer $n$, a connected,
compact, convex, hyperbolic $3$-manifold $C_n(J_{i,j}^-)$ with a local
 isometry (still denoted as
$g_{i,j}$) into $Y_{i,j}$ such  that $C_n(J_{i,j}^-)$
 contains $J_{i,j}^-$ as a hyperbolic submanifold.
The manifold $C_n(J_{i,j}^-)$ is obtained by gluing together $J_{i,j}^-$
 with $n_{i,j}$
``multi-$1$-handles'' $H_{i,j,a}(n), a=1,...,n_{i,j}$, along the attaching region $\p_p J_{i,j}^-$,
where $n_{i,j}$ is the number of components of  $\p S_{i,j}^-$. But there is a subtle difference from
the construction of \cite[Section 8]{MZ} in choosing ``the wrapping numbers'' of the handles
$H_{i,j,a}(n)$.

\begin{adj}\label{diff lengths}
{\rm If $\b$ is a component of $\p S_{i,j}^-$ which lies in the component $T_k$ of $\p M$, then the
multi-$1$-handle associated to it, say the $a$-th one $H_{i,j,a}(n)$,
 will have ``wrapping number''
$\frac{n d_i}{d_{i,k}}$ (instead of $n$ given  in \cite[Section 8]{MZ}),
where $d_{i,k}$ and $d_i$ were
defined  in Notation \ref{dik}.}
\end{adj}

\section{Finding the right covers}\label{lifting}

Recall  the definitions of $n_i$ and $d_i$ given in Notation \ref{dik}.
The main task of this section is
to prove the following

\begin{thm}\label{each large n}Given $S_{i,j}^-$, there is a positive
 even integer $N_{i,j}$ such that for each even integer $N_*\geq N_{i,j}$, we have
\newline (1)  $S_{i,j}^-$ has an
$$m_i=N_* d_i+1$$
fold cover $\breve  S_{i,j}^-$ with $|\p \breve  S_{i,j}^-|=|\p S_{i,j}^-|$ (i.e. each component of $\p
\breve  S_{i,j}^-$ is an $m_i$-fold cyclic cover of a component of $\p  S_{i,j}^-$). So equivalently
each $Y_{i,j}^-$ has an
$$m_i=N_* d_i+1$$
fold cover $\breve  Y_{i,j}^-$ with $|\p_p \breve  Y_{i,j}^-|=|\p_p Y_{i,j}^-|$
(i.e. each component of $\p_p \breve  Y_{i,j}^-$ is an $m_i$-fold cyclic cover of
a component of $\p_p  Y_{i,j}^-$).
\newline
(2) The  map  $g_{i,j}:J_{i,j}^-\ra Y_{i,j}^-$ lifts to an embedding
$\breve g_{i,j}:J_{i,j}^-\ra \breve
Y_{i,j}^-$ and if $\tilde A$ is a component of $\p_p \breve  Y_{i,j}^- $, then components of $\breve
g_{i,j}(\p_p J_{i,j}^-)\cap \breve  A$ are evenly spaced along $\breve  A$. More precisely if $\breve
\b$ is the component of $\p \breve  S_{i,j}^-$
 corresponding to  $\breve  A$,
covering a component $\b$ of $\p S_{i,j}^-$ in $T_k$, then the
topological center points of $\tilde
g_{i,j}(\p_p J_{i,j}^-)\cap \breve  A$ divide $\breve  \b$ into arc
components each with wrapping number
$N_*\frac{d_i}{d_{i,k}}$.
\end{thm}

Of course in Theorem \ref{each large n}, the cover $\tilde S_{i,j}^-$ and the number
$m_i$ depend on $N_*$.
For  simplicity, we suppressed this dependence
in notation for $\tilde S_{i,j}^-$ and $m_i$.
Similar  suppressed
notations will occur also in other places  later in the paper when there is no
danger of causing confusion, and we shall not remark on this
all the time.

For the definition of the wrapping number see Definition \ref{wrapping}.

\begin{cor}\label{indep of j}There is a positive even integer $N_0$  such that
for each even integer $N_*\geq N_0$  and for  each $i=1,2, j=1,...,n_i$,
we have
 \newline
 (1)  $S_{i,j}^-$ has an
$$m_i=N_* d_i+1$$
fold cover $\breve  S_{i,j}^-$ with $|\p \breve  S_{i,j}^-|=|\p S_{i,j}^-|$.
 So equivalently each $Y_{i,j}^-$ has an
$$m_i=N_* d_i+1$$
fold cover $\breve  Y_{i,j}^-$ with $|\p_p \breve  Y_{i,j}^-|=|\p_p Y_{i,j}^-|$.
\newline
(2) The  map  $g_{i,j}:J_{i,j}^-\ra Y_{i,j}^-$ lifts to an embedding $\breve g_{i,j}:J_{i,j}^-\ra \breve
Y_{i,j}^-$ and if $\tilde A$ is a component of $\p_p \breve  Y_{i,j}^- $, then components of $\breve
g_{i,j}(\p_p J_{i,j}^-)\cap \breve  A$ are evenly spaced along $\breve  A$. More precisely if $\breve
\b$ is the component of $\p \breve  S_{i,j}^-$
 corresponding to  $\breve  A$,
covering a component $\b$ of $\p S_{i,j}^-$ in $T_k$, then the topological center points of $\tilde
g_{i,j}(\p_p J_{i,j}^-)\cap \breve  A$ divide $\breve  \b$ into arc components each with wrapping number
$N_*\frac{d_i}{d_{i,k}}$.
\end{cor}

\pf Apply Theorem \ref{each large n} and let $N_0=max\{N_{i,j}; i=1,2,j=1,...,n_i\}$.\qed

Corollary \ref{indep of j} is to say that the number $N_*$ and thus the number $m_i$ are independent of
the second index $j$ in $S_{i,j}^-$.

For notational simplicity, we shall only consider the following two cases
in proving Theorem \ref{each large n}:

\noindent{\bf Case 1}. A given surface $S_{i,j}^-$ has $b_1$
boundary components $\{\b_{1,p}, p=1,...,
b_1\}$ on $T_1$ and $b_2$ boundary components $\{\b_{2,p}, p=1,...,b_2\}$ on $T_2$, and is disjoint from
$T_3,...,T_m$. So $n_{i,j}=b_1+b_2$ which is the number of components
 of $\p S_{i,j}^-$.

\noindent{\bf Case 2}. A given surface $S_{i,j}^-$ has only one  boundary component  $\{\b\}$ on $T_1$
and is disjoint from  $T_2,...,T_m$. So $n_{i,j}=1$, which is the number of components of $\p
S_{i,j}^-$.

The reader will see that the proof of Theorem \ref{each large n} for a general surface $S_{i,j}^-$
 will be very similar to either case 1 or 2, depending on whether $S_{i,j}^-$ has multiple boundary components,
 or just a single one.

 \noindent{\bf Proof of
Theorem \ref{each large n} in Case 1}.

Again to avoid too complicated notations on indices,
 in the following we shall suppress the indices
$i,j$
 for some items  depending on them, when there is no danger of
 causing confusion.

 Recall that $\p S_{i,j}^-$
have the  induced orientation from  the orientation of $S_{i,j}^-$.
Let $\b_{k,p}, p=1,...,b_k$ be the oriented  boundary components of $\p S_{i,j}$
in $T_k$ for each $k=1,2$.
Recall the number $d_{i,k}$ given in Notation \ref{dik}.
We list the set of intersection points
of $\p S_{i,j}^-$ with $\p S_{i_*}^-$  as $t_{k,p,q}$,
$k=1,2$, $p=1,...,b_k$ and $q=1,...,d_{i,k}$, so that $t_{k,p,q}, q=1,...,d_{i,k}$, appear consecutively
along $\b_{k,p}$ following its orientation.
 We choose $t_{1,1,1}$
as the base point for $\pi_1(S_{i,j}^-)=\pi_1(Y_{i,j}^-)=\pi_1(S_{i,j})=\pi_1(Y_{i,j})$.

\begin{note}\label{n even}
{\rm {\bf From now on in this paper we assume that $n$ is a positive  even number}}
\end{note}

 Recall that
there is a local isometry $g_{i,j}:C_n(J_{i,j}^-)\ra Y_{i,j}$ which is a one-to-one map  when restricted
to the set of center points of $\p_p J_{i,j}^-$. We list these center points as $b_{k,p,q}$ so that
$t_{k,p,q}=g_{i,j}(b_{k,p,q})$ for all $k,p,q$. We choose $b_{1,1,1}$ as the  base point for each of
$J_{i,j}$, $J_{i,j}^-$, $\hat J_{i,j}$ and $C_n(J_{i,j}^-)$.

Similar to  the definition given in \cite[p2144]{MZ}, we have

\begin{defi}\label{wrapping}{\rm  Suppose that
 $\breve p:\breve{\b}_{k,p}\ra
\b_{k,p}$ is a covering
 map, and let $\breve \b_{k,p}$ have the orientation induced from that of $\b_{k,p}$.
 Let $\a\subset \tilde \b_{k,p}$ be an embedded, connected, compact arc
 with the orientation induced from that of $\breve \b_{k,p}$, whose
 initial point is in $\breve p^{-1}(t_{k,p,q})$
 and whose terminal point is in
 $\breve p^{-1}(t_{k,p,q+1})$ (here $q+1$ is defined mod $d_{i,k}$).
 We say that $\a$ has {\it wrapping number} $n$ if
 there are exactly
$n$ distinct points of $\breve p_i^{-1}(t_{k,p,q})$ which are contained in the
interior of $\a$.}
\end{defi}

\begin{nota}
\label{generators}
{\rm Let $g$ be the genus of $S_{i,j}^-$. As in \cite[Section 10]{MZ}, the group $\pi_1(S_{i,j}^-, t_{1,1,1})$ has a set
of generators
$$L=\{a_1,b_1,a_2, b_2, ...,a_g,b_g,
x_{1},x_2, ...,x_{n_{i,j}-1}\}$$ such that the elements
$$x_{1},x_{2},...,x_{n_{i,j}-1},x_{n_{i,j}}=
[a_{1},b_{1}][a_{2},b_{2}]\cdots [a_{g},b_{g}]x_{1}x_{2}\cdots x_{n_{i,j}-1}$$
 have
representative loops, based at the point $t_{1,1,1}$,   freely homotopic to  the
$n_{i,j}=b_1+b_2$  components
$\b_{1,1},\b_{1,2},...,\b_{1,b_1},\b_{2,1},\b_{2,2},...,\b_{2,b_2}$ of $\p
S_{i,j}^-$ respectively.
}
\end{nota}

As in \cite[Section 10]{MZ}, we fix  a  generating set  $$w_1,...,w_\ell$$ for $\pi_1(J_{i,j}^-,
b_{1,1,1})$ and choose a generating set
 $$w_1,...,w_\ell, z_{k,p,q}(n), \;k=1,2, p=1,..,b_k, q=1,...,d_{i,k}-1$$
 for $\pi_1(C_n(J_{i,j}^-),b_{1,1,1})$ such that
$$\pi_1(C_n(J_{i,j}^-),b_{1,1,1})=\pi_1(J_{i,j}^-, b_{1,1,1})*<z_{k,p,q}(n),
\;k=1,2, p=1,..,b_k, q=1,...,d_{i,k}-1>$$ where $*$ denotes the free product, and $<z_{k,p,q}(n),
\;k=1,2, p=1,..,b_k, q=1,...,d_{i,k}-1>$ is the free group freely generated by the $z_{k,p,q}(n)$'s.

Here are some necessary details of how $z_{k,p,q}(n)$ is defined, following \cite[Section 10]{MZ} but
with different and simplified notations for indices.
 Let $\alpha_{k,p,q}\subset J_{i,j}^-$ be a fixed,
 oriented path from $b_{1,1,1}$ to $b_{k,p,q}$,
 for each of $k=1,2, p=1,...,b_k, q=1,...,d_{i,k}$ ($\a_{1,1,1}$ is the constant path).
 Recall the construction of $C_n(J_{i,j}^-)$
and Adjustment \ref{diff lengths}. For $k=1,2$, $p=1,...,b_k$,   let $H_{k,p}(n)$ denote the
multi-1-handle of $C_n(J_{i,j}^-)$ corresponding to the component $\b_{k,p}$ of $\p S_{i,j}^-$.
 For $k=1,2$, $p=1,..., b_k$, $1 \leq q \leq d_{i,k}-1$,
 let $\delta_{k,p,q}(n)$ be the oriented geodesic arc in the multi-one-handle
  $H_{k,p}(n)\subset C_n(J_{i,j}^-)$ from the point
 $b_{k,p,q}$ to $b_{k,p,q+1}$.
 Then $$z_{k,p,q}(n)=
 \alpha_{k,p,q}\cdot\delta_{k,p,q}(n)\cdot\overline{\alpha_{k,p,q+1}},$$
 where the symbol ``$\cdot$'' denotes path concatenation
 (sometimes omitted), and $\overline{\alpha_{k,p,q+1}}$
 denotes the reverse of $\alpha_{k,p,q+1}$.
 Also we always write path (in particular loop) concatenation from left to right.

As in \cite[Section 10]{MZ}, if $\a$ is an oriented arc in $C_n(J_{i,j}^-)$, we use $\a^*$ to denote the
oriented arc $g_{i,j}\circ\a$ in $Y_{i,j}$, and if $\g$ is an element in $\pi_1(C_n(J_{i,j}^-,
b_{1,1,1})$, we use $\g^*$ to denote the element $g_{i,j}^*(\g)$ where $g_{i,j}^*$ is the induced
homomorphism $g_{i,j}^*: \pi_1( C_N(J_{i,j}^-), b_{1,1,1})\ra \pi_1(Y_{i,j}, t_{1,1,1})$.

The oriented path $\alpha_{k,p,q}^*$ in $Y_{i,j}^-$ runs from $t_{1,1,1}$ to
$t_{k,p,q}$.
 For $k=1,2$, $p=1,...,b_k$, $1 \leq q \leq d_{i,k}-1$,
 let $\eta_{k,p,q}$ be the oriented subarc in $\beta_{k,p}$ from
 $t_{k,p,q}$ to $t_{k,p,q+1}$ following the orientation of $\b_{k,p}$, and
 let $\s_{k,p,q}\subset Y_{i,j}^-$ be the loop
 $\alpha_{k,p,q}^*\cdot\eta_{k,p,q}\cdot\overline{\alpha_{k,p,q+1}}^*$.
  Let $\s_{k,p,0}$ be the constant path based at $t_{1,1,1}$.
 Let $x_{b}^{\prime}$ be the loop
 $\alpha_{1,b,1}^*\cdot\beta_{1,b}\cdot\overline{\alpha_{1,b,1}}^*$
 if $b=1,...,b_1$ and be the loop
  $\alpha_{2,b-b_1,1}^*\cdot\beta_{2,b-b_1}\cdot\overline{\alpha_{2,b-b_1,1}}^*$
 if $b=b_1+1,..., b_1+b_2=n_{i,j}$,
 where  $\b_{k,p}$ is considered as an oriented loop
 starting and ending at the point $t_{k,p,1}$.
Similar to \cite[Lemma 10.1]{MZ}, we have

\begin{lem}\label{diff powers}
\label{lk} Considered as an element in $\pi_1(Y_{i,j}, t_{1,1,1})=\pi_1(S_{i,j}, t_{1,1,1})
=\pi_1(S_{i,j}^-, t_{1,1,1})$,
  $$z_{k,p,q}(n)^* =
  (\overline{\s_{k,p,q-1}}\cdots\overline{\s_{k,p,0}})
  (x_{b_{k-1}+p}')^{n\frac{d_i}{d_{i,k}}}
 (\s_{k,p,0}\cdots\s_{k,p,q}),
$$ for each of $k=1,2, p=1,...,b_k, q=1,...,d_{i,k}-1$, where $b_0$ is defined to be $0$.\qed
\end{lem}

\begin{note}{\rm The only real meaningful difference of this lemma from \cite[Lemma 10.1]{MZ} is that the power
of $x_{b_{k-1}+p}'$ in the expression of $z_{k,p,q}(n)^*$ above depends  on the indices  $i$ and $k$ besides
$n$, which is due to the Adjustment \ref{diff lengths}. }
\end{note}

Recall  that $\hat J_{i,j}$ is a connected, compact, convex, hyperbolic
3-manifold obtained from $J_{i,j}^-$ by capping off each component of $\p_p
J_{i,j}^-$ with a compact, convex $3$-ball, and that $\pi_1(J_{i,j},
b_{1,1,1})=\pi_1(J_{i,j}^-,b_{1,1,1})=\pi_1(\hat J_{i,j},b_{1,1,1})$.
 Also, $\hat J_{i,j}$ is a
submanifold of $C_n(J_{i,j}^-)$, so by \cite[Lemma 4.2]{MZ}, $\pi_1(\hat J_{i,j},b_{1,1,1})$ can be
considered as a subgroup of $\pi_1(C_n(J_{i,j}^-),b_{1,1,1})$. As
 $C_n(J_{i,j}^-)$ is a connected, compact, convex, hyperbolic
3-manifold, the induced homomorphism $g_{i,j}^*:\pi_1(C_n(J_{i,j}^-), b_{1,1,1})\ra \pi_1(Y_{i,j},
t_{1,1,1})
 =\pi_1(S_{i,j}^-, t_{1,1,1})$ is
 injective by again \cite[Lemma 4.2]{MZ}.
So  $g_{i,j}^*(\pi_1(C_n(J_{i,j}^-), b_{1,1,1}))= g_{i,j}^*(\pi_1(J_{i,j}^-, b_{1,1,1}))
*<z_{k,p,q}(n),
\;k=1,2, p=1,..,b_k, q=1,...,d_{i,k}-1>$ is a subgroup of $\pi_1(Y_{i,j}, t_{1,1,1})=\pi_1(S_{i,j},
t_{1,1,1}) =\pi_1(S_{i,j}^-, t_{1,1,1})$.

By \cite[Proposition 4.7]{MZ}
 there is a set of elements $y_1,..., y_r$ in
$$\pi_1(Y_{i,j}, t_{1,1,1}) -g_{i,j}^*(\pi_1(\hat J_{i,j},b_{1,1,1}))$$ such that,
if $G$ is a finite
index subgroup of $\pi_1(Y_{i,j}, t_{1,1,1})$ which separates
$g_{i,j}^*(\pi_1(\hat J_{i,j},b_{1,1,1}))$
from $y_1,...,y_r$, then the local isometry $g_{i,j}:\hat J_{i,j}\ra Y_{i,j}$
lifts to an embedding
$\breve g_{i,j}$ in the finite cover $\breve Y_{i,j}$ of $Y_{i,j}$
corresponding to $G$.

To prove Theorem \ref{each large n} in Case 1,  we just need to prove the following

\begin{thm}\label{large n}
There is a positive even integer $N_{i,j}$ such that for each  even integer $N_*\geq N_{i,j}$,
there is a finite index
subgroup $G$ of $\pi_1(Y_{i,j}, t_{1,1,1})=\pi_1(S_{i,j}, t_{1,1,1}) =\pi_1(S_{i,j}^-, t_{1,1,1})$
 such
that
\newline (i) $G$ has index $m_i=N_*d_i+1$;
\newline (ii) $G$ contains the elements
$w_{1}^*,...,w_{l}^*$, and thus contains the subgroup $g_{i,j}^*(\pi_1(\hat
J_{i,j},b_{1,1,1}))=g_{i,j}^*(\pi_1(\hat J_{i,j}^-,b_{1,1,1}))$;
\newline
(iii) $G$ contains the elements $z_{k,p,q}(N_*)^*$, $k=1,2$, $p=1,...,b_k$, $q=1,...,d_{i,k}-1$;
\newline (iv) $G$ does not contain any of
$x_b^d$, $b=1,...,n_{i,j}$, and $d=1,...,m_i-1$;
\newline
(v) $G$ does not contain any of $y_1,...,y_r$.
\end{thm}

\begin{pro}
 Theorem \ref{each large n} in Case 1 follows from
 Theorem \ref{large n}.
 \end{pro}

 \pf The  proof is similar to that of \cite[Proposition 11.1]{MZ}.

Let $\breve Y_{i,j}$ be the finite cover of $Y_{i,j}$ corresponding the subgroup $G$ provided by
Theorem
\ref{large n}, and let $\breve S_{i,j}$ be the corresponding center surface of
 $\breve Y_{i,j}$ covering $S_{i,j}$.

As noted earlier, Conditions (ii) and (v) of Theorem \ref{large n}
 imply that the map $g_{i,j}:\hat J_{i,j}\ra Y_{i,j}$ lifts
 to an embedding $\breve g_{i,j}:\hat J_{i,j}\ra \breve Y_{i,j}$.
 Conditions (i) and (iv) of
Theorem \ref{large n} imply  that $|\p \breve S_{i,j}^-|=|\p S_{i,j}|$. So part (1) of Theorem \ref{each
large n} holds in Case 1.

We may now let $\breve \b_{k,p}$ be the component of $\p \breve S_{i,j}^-$ covering $\b_{k,p}$ for each
of $k=1,2, p=1,..., b_k$.

 Conditions (ii) and (iii) of
Theorem \ref{large n} imply that
 the group $g_{i,j}^*(\pi_1(C_{N^*}(J_{i,j}^-), b_{1,1,1}))$ is
contained in $G$.
 Therefore the map
 $g_{i,j}:(C_{N_*}(J_{i,j}^-),b_{1,1,1})\ra (Y_{i,j}, t_{1,1,1})$ lifts
 to a map
 $$\breve g_{i,j}:(C_{N_*}(J_{i,j}^-), b_{1,1,1})\ra
 (\breve Y_{i,j}, \breve g_{i,j}(b_{1,1,1})).$$
 Recall the notations established earlier.
 Consider the multi-1-handle $H_{k,p}(N_*) \subset C_{N_*}(J_{i,j}^-)$
 containing the points $b_{k,p,q}, q=1,...,d_{i,k}$, and the
 geodesic arcs $\d_{k,p, q}(N_*)\subset H_{k,p}(N_*)$, $q=1,...,d_{i,k}-1$.
 By our construction the immersed arc $g_{i,j}:\d_{k,p,q}(N_*)\ra S_{i,j}$
 is homotopic, with end points fixed, to the path
 in $\b_{k,p}$ which starts at the point $t_{k,p,q}$,
 wraps $N_*\frac{d_i}{d_{i,k}}$ times around $\b_{k,p}$
  and then continues to the point
 $t_{k,p,q+1}$, following the orientation of $\b_{k,p}$.
 This latter (immersed) path lifts to an embedded arc in
$\breve \b_{k,p}$ connecting $\breve g_{i, j}(b_{k,p,q})$ and $\breve g_{i, j}(b_{k,p,q+1})$, because
$\breve \b_{k,p}$ is an $N_*d_i+1$-fold cyclic cover of $\b_{k,p}$. This shows that part (2) of Theorem
\ref{each large n} holds in Case 1. \qed

Theorem \ref{large n} is proved using the technique of folded graphs. We shall follow as closely as
possible the approach used in \cite{MZ} and we assume the terminologies
 used there
concerning $L$-directed graphs. Recall that $L$ is the generating set
chosen in Notation \ref{generators} for the free group $\pi_1(Y_{i,j},
t_{1,1,1})=\pi_1(S_{i,j}, t_{1,1,1}) =\pi_1(S_{i,j}^-, t_{1,1,1})$. From now on any
 group element in  $\pi_1(Y_{i,j},
t_{1,1,1})=\pi_1(S_{i,j}, t_{1,1,1}) =\pi_1(S_{i,j}^-, t_{1,1,1})$ will be considered as a word in
$L\cup L^{-1}$.

First we translate Theorem \ref{large n} into the following theorem, in terms of folded graphs.

\begin{thm}\label{large} There is a positive even integer $N_{i,j}$ such that for
each even integer $N_*\geq N_{i,j}$
there is a finite, connected, $L$-labeled, directed graph ${\cal G}(N_*)$
 (with a fixed base vertex
$v_0$) with the following properties:
\newline
  (0)  ${\cal G}(N_*)$ is  $L$-regular;\newline
(1) The number of vertices of  ${\cal G}(N_*)$ is
 $m_{i}=N_*d_i+1$; \newline
 (2) Each of the words $w_1^*,..., w_\ell^*$ is  representable
 by a loop, based at $v_0$, in ${\cal G}(N_*)$; \newline
(3) ${\cal G}(N_*)$  contains a closed loop, based at $v_0$,
representing the word $z_{k,p,q}(N_*)^*$,
for each $k=1,2, p=1,...,b_k, q=1,...,d_{i,k}-1$;\newline
 (4)  ${\cal G}(N_*)$ contains no
 loop representing the word $x_b^d$ for any $b=1,...,n_{i,j}$ and
 $d=1,...,m_{i}-1$;\newline
  (5) each of the words $y_1,...,y_r$ is  representable by a
non-closed path, based at  $v_0$, in ${\cal G}(N_*)$.
\end{thm}

\begin{note}{\rm  The procedure for constructing the graphs
described in Theorem \ref{large} follows mostly  that given in \cite[Section 11]{MZ}.
In the current case we need to deal with two  major complications.
 One is due to the fact that the number $d_{i, k}$  of intersection points
 in  a boundary   component $\b_{k,p}$ of $\p S_{i,j}^-$  depends on $k$;
 the other is due to the requirement of showing that such a graph ${\cal G}(N_*)$ exists
  for each even integer $N_*\geq N_{i,j}$. Actually our adjustment has begun as early as in
 Adjustment  \ref{diff lengths}.
 }
\end{note}

\cite[Remark 9.7]{MZ} was one of the main group theoretical results  obtained
in \cite{MZ} and it will
also play a fundamental role in our current case. We quote this result below
as Theorem \ref{regular} in
the current notations.

\begin{thm}\label{regular}{\rm (\cite[Remark 9.7]{MZ})}
If ${\cal G}_\#$ is a finite, connected, $L$-labeled, directed,
folded graph with base vertex $v_0$,
with corresponding
 subgroup $G_\# = L({\cal G}_\#,v_0) \subset \pi_1(S_{i,j}^-,t_{1,1,1})$,
 such that
\\
$\bullet\;$ ${\cal G}_\#$ does not contain any loop representing the word $x_b^d$ for any
$b=1,...,n_{i,j}$, $d\in \z-\{0\}$, and
\\
$\bullet\;$ $y_1,...,y_r$ are some fixed, non-closed paths
 based at $v_0$ in ${\cal G}_\#$,\\
 then there is a finite, connected, $L$-regular
graph ${\cal G}_*$  such that
\\
$\bullet\;$ ${\cal G}_*$ contains ${\cal G}_\#$ as an embedded subgraph,
 and thus in particular $y_1,...,y_r$ remain non-closed paths
 based at $v_0$ in ${\cal G}_*$, and
\\
$\bullet\;$ ${\cal G}_*$ contains no loops representing the word $x_b^d$, for
each of $b=1,...,n_{i,j}$,
$d=1,...,m_*-1$, where $m_*$ is the number of vertices of ${\cal G}_*$. \qed
\end{thm}

We now begin our constructional proof of Theorem \ref{large}.
Let ${\cal G}_{1}$ be the connected,
finite, $L$-labeled, directed graph which
 results from taking a disjoint union of embedded loops-- representing
the reduced versions of the words $w_{1}^*, ..., w_{\ell}^*$
  respectively-- and non-closed embedded paths--
 representing the reduced versions of the
 words $y_{1}, ..., y_{r}$
  respectively--
 and then identifying their base points to
 a common vertex $v_{0}$.
 Then $L({\cal G}_{1}, v_{0})$ represents the subgroup
$g_{i,j}^*(\pi_1( J_{i,j}^-, b_{1,1,1}))$ of $\pi_1(S_{i,j}^-, t_{1,1,1})$.
Since the folding operation
does not change the group that the graph represents,
 $L({\cal G}_{1}^f, v_{0})=g_{i,j}^*(\pi_1(J_{i,j}^-,b_{1,1,1}))$
 (where ${\cal G}_1^f$ denotes the folded graph of ${\cal G}_1$).

Recall from Lemma \ref{diff powers} that  $$z_{k,p,q}(n)^* =
  (\overline{\s_{k,p,q-1}}\cdots\overline{\s_{k,p,0}})(x_{b_{k-1}+p}')^{n\frac{d_i}{d_{i,k}}}
 (\s_{k,p,0}\cdots\s_{k,p,q}),$$
 $k=1,2, p=1,...,b_k, q=1,...,d_{i,k}-1$.
 Note that $x_{b}'$ is conjugate to $x_b$  in $\pi_1(S_{i,j}^-, t_{1,1,1})$,
   for $b=1,..., n_{i,j}$.
 Let $\t_{b}$ be an element of $\pi_1(S_{i,j}^-, t_{1,1,1})$ such that
 $x_{b}^{\prime} = \t_{b} x_{b}\t_{b}^{-1}$, $b=1,..., n_{i,j}$.
 Let ${\cal G}_{2}$ be the connected graph which
 results from taking the disjoint union of  ${\cal G}_{1}^f$
 and non-closed embedded paths representing  the reduced version of
 the words  $\overline{\s_{k,p,q-1}}\cdots\overline{\s_{k,p,0}}\t_{b_{k-1}+p}$,
  $k=1,2, p=1,...,b_k, q=1,...,d_{i,k}-1$, respectively, and then identifying their base
vertices into a single base vertex which we still denote by
$v_{0}$. Then obviously we have $L({\cal G}_{2}^f, v_{0}) =
L({\cal G}_{2}, v_{0})
    = L({\cal G}_{1}^f, v_{0})=g_{i,j}^*(\pi_1(J_{i,j}^-,b_{1,1,1}))$.

 Let
 $v_{k,p,q}$ be the terminal vertex
 of the path  $\overline{\s_{k,p,q-1}}\cdots\overline{\s_{k,p,0}}\t_{b_{k-1}+p}$
 in ${\cal G}_{2}^f$, for each $k=1,2, p=1,...,b_k, q=1,...,d_{i,k}$.
For each  of $(1, p, q)$, $p=1,...,b_1, q=1,...,d_{i,1}$, and $(2, p,q)$, $p=1,..., b_2-1$,
 $q=1,...,d_{i,2}$, let $Q_{k,p,q}$ be the maximal
 $x_{b_{k-1}+p}$-path in $\widehat{{\cal G}_{2}^f}$ (a maximal $x_{b}$-path was
 defined in \cite[Section 9]{MZ}) which contains the vertex $v_{k,p,q}$.
 For each of $(2, b_2, q)$, $q=1,..., d_{i,2}$, let $Q_{2,b_2,q}$  be the
 maximal $x_{n_{i,j}}$-path  in $\widehat{{\cal G}_{2}^f}$
determined by \newline (1) if there is a directed edge of
$\widehat{{\cal G}_{2}^f}$ with
 $v_{2,b_2,q}$ as its initial vertex
 and with the first letter of the word $x_{n_{i,j}}$ as its label, then $Q_{2, b_2,q}$
 contains that edge;
\newline
(2) if the edge described in (1) does not exists, then $v_{2,b_2,q}$
is the terminal vertex of  $Q_{2,b_2,q}$ and  the first letter of
the word $x_{n_{i,j}}$ is the
 terminal missing label of $Q_{2,b_2,q}$.
\newline
 Note that each $Q_{k, p,q}$ is uniquely determined.
 Also no  $Q_{k, p,q}$ can
 be an $x_{b}$-loop, since the group
 $L({\cal G}_{2}^f, v_{0})=g_{i,j}^*(\pi_1(J_{i,j}^-,b_{1,1,1}))$
      does not contain non-trivial peripheral elements of $\pi_1(S_{i,j}^-, t_{1,1,1})$
      by \cite[Lemma 4.2]{MZ}.
  Let  $v_{k,p,q}^-$ and $v_{k,p,q}^+$ be the initial
and terminal
 vertices of $Q_{k,p,q}$ respectively.
Note that if $p\ne b_2$ and  $Q_{k,p,q}$ is not a constant
 path, then $v_{k,p,q}^-$ and $v_{k,p,q}^+$ must be
distinct vertices; however
  $v_{2,b_2, q}^-$ and $v_{2,b_2,q}^+$ may possibly be the
 same vertex, even if $Q_{2, b_2, q}$ is a non-constant path.

Let $Q_{k,p,q}^-$ be the embedded subpath of $Q_{k,p,q}$ with
$v_{k,p,q}^-$ as the initial vertex and with $v_{k,p,q}$ as the
terminal vertex, and  let $Q_{k,p,q}^+$ be the embedded subpath of
$Q_{k,p,q}$ with $v_{k,p,q}$ as the initial vertex and with
$v_{k,p,q}^+$ as the terminal vertex.

 Note
that the number
 $max\{Length(Q_{k,p,q}): k=1,2, p=1,...,b_k, q=1,...,d_{i,k}\}$ is independent of $n$,
and thus is bounded. So we may assume that
$$n>40|L|+2 max\{Length(Q_{k,p,q}): k=1,2, p=1,...,b_k, q=1,...,d_{i,k}\}.$$
We shall also assume that $n$ has been chosen large enough so that
$C_n(J_{i,j}^-)$ is convex.

 Now for each  $k=1,2, p=1,...,b_k$, $q=1,...,d_{i,k}-1$,  we make a new
 non-closed embedded path $\Theta_{k,p,q}(n)$
representing the word $x_{b_{k-1}+p}^{n\frac{d_i}{d_{i,k}}}$, and we add it to the graph ${\cal
G}_{2}^f$, by identifying the initial vertex of $\Theta_{k,p,q}(n)$ with $v_{k,p,q}$
 and the terminal
vertex
  with $v_{k,p,q+1}$.

\begin{adj}{\rm For each  $k=1,2, p=1,...,b_k$,
$q=d_{i,k}$, we make a new non-closed embedded path $\Theta_{k,p,q}(n)$ representing the word
$x_{b_{k-1}+p}^{n\frac{d_i}{d_{i,k}}}$, and we add it to the graph ${\cal G}_{2}^f$, by identifying the
initial vertex of $\Theta_{k,p,q}(n)$ with $v_{k,p,q}$. }
\end{adj}

\begin{note}\label{total}{\rm  For each fixed $k=1,2, p=1,...,b_k$, the paths
$\{\Theta_{k,p,q}(n), q=1,...,d_{i,k}\}$, are connected together and
 form a connected path representing the word $x_{b_{k-1}+p}^{n d_i}$.}
\end{note}

 In the resulting graph there are some obvious
 places one can perform the folding operation:
 for each  $k=1,2, p=1,...,b_k, q=1,...,d_{i,k}-1$, the path
$Q_{k,p,q}^+$ can be completely folded into the added new path $\Theta_{k,p,q}(n)$, and likewise the
path $Q_{k,p,q+1}^-$ can be completely folded into $\Theta_{k,p,q}(n)$.
 Let ${\cal G}_{3}(n)$ be the
resulting graph after performing these specific
 folding operations  for each  $k=1,2, p=1,...,b_k, q=1,...,d_{i,k}$.

 From the explicit construction, it is clear that ${\cal G}_{3}(n)$
 has the following properties:
 \newline
 (1) ${\cal G}_{3}(n)$ is a connected, finite, $L$-labeled,
 directed graph;
\newline
(2) ${\cal G}_{3}(n)$ contains loops, based at $v_{0}$, representing the word $z_{k,p,q}(n)^*$ for each
$k=1,2, p=1,...,b_k, q =1,...,d_{i,k}-1$;
 \newline
 (3) ${\cal G}_{3}(n)$  contains ${\cal G}_{2}^f$ as an embedded
 subgraph.

It follows from Property (3) that the paths in ${\cal G}_{2}^f$ representing the
words $y_{1}, ..., y_{r}$ remain each non-closed
  in ${\cal G}_{3}(n)$, and it follows from
 Properties (2) and  (3) and the construction that $L({\cal G}_{3}(n),v_{0})
 =g_{i,j}^*(\pi_1(C_n(J_{i,j}^-), b_{1,1,1}))$. So
$\widehat{{\cal G}_{3}(n)}$ cannot have $x_{b}$-loops for any $b=1,...,n_{i,j}$ (again by \cite[Lemma
4.2]{MZ}).

Now we consider the remaining folding operations on ${\cal G}_{3}(n)$ that need to be done, in order to
get the folded graph ${\cal G}_{3}(n)^f$.

For each  $k=1,2, p=1,...,b_k, q=1,...,d_{i,k}-1$, let
$$\Theta_{k,p,q}(n)' =\Theta_{k,p,q}(n)\setminus
(Q_{k,p,q}^+\cup Q_{k,p,q+1}^-),$$ and for each $k=1,2, p=1,...,b_k, q=d_{i,k}$, let
$$\Theta_{k,p,q}(n)' =\Theta_{k,p,q}(n)\setminus Q_{k,p,q}^+.$$
Then by our construction  each
 $\Theta_{k,p,q}(n)'$ is an embedded
$x_{b_{k-1}+p}$-path with $v_{k,p,q}^+$ as its initial vertex and with $v_{k,p,q+1}^-$ (when $q\ne
d_{i,k}$) as the terminal vertex. Also all these paths $\Theta_{k,p,q}(n)'$,$k=1,2, p=1,...,b_k,
q=1,...,d_{i,k}$, are mutually disjoint in their interior, and their disjoint union is equal to ${\cal
G}_{3}(n)\setminus {\cal G}_{2}^f$.
 Since
$\widehat{{\cal G}_{3}(n)}$ has no $x_{b}$-loops, we see immediately that when  $(k,p)\ne (2, b_2)$, all
the vertices
 $v_{k,p,q}^{\pm}$, $q=1,...,d_{i,k}$, are mutually distinct.

 It follows that the only remaining folds are at the vertices
 $v_{2, b_2,q}^{\pm}$.  At such a
vertex there is at most one edge from $\Theta_{2,b_2,q}(n)'$ which may be folded with one
$x_{b_{k-1}+p}$-edge of $\Theta_{k,p,q_*}(n)'$ at its initial or terminal vertex,
 for some $(k,p)\ne (2,
b_2)$ and some $1\leq q_*\leq d_{i,k}$.
 Thus ${\cal G}_{3}(n)^f$
 is obtained from ${\cal G}_{3}(n)$ by performing at most
 $2d_{i,2}$ folds  (which occur at some of the vertices
$v_{2, b_2,q}^{\pm}$, $q=1,...,d_{i,2}$),
  and every non-closed, reduced path
 in  ${\cal G}_{3}(n)$ which is based at $v_{0}$ will
 remain non-closed in ${\cal G}_{3}(n)^f$.
 In particular, the paths representing the words $y_{1}, ... y_{r}$ are
 each non-closed in
 ${\cal G}_{3}(n)^f$.

Let $f_3:{\cal G}_{3}(n)\ra {\cal G}_{3}(n)^f$ be the natural map and we fix a number
 $$s>
 2d_{i,2}+Diameter({\cal G}_{2}^f).$$ Then the  map
 $f_3: {\cal G}_{3}(n) \ra {\cal G}_{3}(n)^f$
 is an embedding on ${\cal G}_{3}(n) - N_{s}(v_{0})$,
 where $N_s(v_0)$ denotes the $s$-neighborhood of $v_0$ in ${\cal G}_3(n)$
  considering a graph  as
a metric space, by making each
 edge isometric to the interval $[0,1]$.
 Obviously the number $s$ is independent of $n$.

 \begin{note}\label{large power} {\rm
    We may assume further that  $n$ is large enough
 so that
 the components of ${\cal G}_{3}(n)^f\setminus f(N_{v_{0}}(s))$
can be denoted by   $\Phi_{k,p,q}(n)$, $k=1,2, p=1,...,b_k, q=1,...,d_{i,k}$,
 such that  $\Phi_{k,p,q}(n)$ is an embedded subpath in
 $\Theta_{k,p,q}(n)'$  containing
 a power of $x_{b_{k-1}+p}$ which is larger than
 $\frac{n d_i }{d_{i,k}}-\frac{n}{4}$.
 This is clearly possible from the construction.}
\end{note}

 The next step is to modify the graph  ${\cal G}_{3}(n)^f$,
 by inserting copies of a certain graph $\Omega$,
 pictured in Figure \ref{fg1}, and then performing folding
 operations, to obtain a  graph  (the graph ${\cal
  G}_{4}(n)$ given below)
  which contains loops, based at the base vertex $v_{0}$,
  representing the words  $$w_{1}^*, ...,
w_{\ell}^*, z_{k,p,q}(n+1)^*, k=1,2, p=1,...,b_k,q=1,..., d_{i,k}-1,$$ respectively,
 and which contains
 non-closed paths, based at $v_{0}$, representing the
 words $$y_{1}, ..., y_{r}$$ respectively.
  In Figure \ref{fg1}, single edge loops at a vertex have label
one
 each from $L^*=\{a_{1}, b_{1},..., a_{g}, b_{g}\}$.
The edges in part (a) and (b)  connecting two adjacent vertices are $x_{b}$-edges,
$b=1,2,...,n_{i,j}-1$, (precisely $n_{i,j}-1$ edges). In part (a) of the figure, an $x_{b}$-edge points
from the left vertex to the right vertex iff
 $b$ is odd, and in part (b) of the figure,
 an $x_{b}$-edge points from  left to
right iff
 $b$ is $1$ or an even number.

\begin{figure}[!ht]
{\epsfxsize=4.5in \centerline{\epsfbox{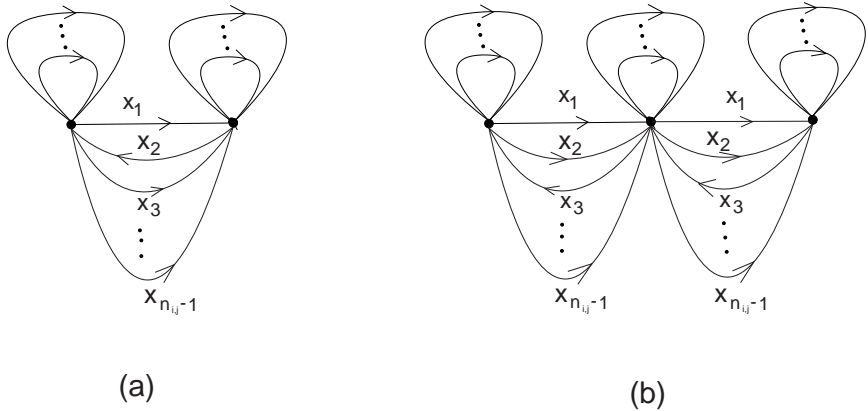}}\hspace{10mm}} \caption{The graph $\Omega$ when
 (a) $n_{i,j}>0$ is even, (b) $n_{i,j}>1$ is odd.}\label{fg1}
\end{figure}

 The method for constructing ${\cal G}_{4}(n)$
  breaks into two subcases, i.e.
 \newline(a) when $n_{i,j}$ is even, \newline(b) when $n_{i,j}$ is odd.

 \noindent
 {\bf Subcase (a):}  $n_{i,j}$ is even.

 Recall that for each $k=1,2$, $p=1,...,b_k$,
 $\cup_{q=1}^{d_{i,k}}\Theta_{k,p,q}(n)$ is a connected path in
 ${\cal{G}}_{3}(n)^f$ representing the word
 $$x_{b_{k-1}+p}^{nd_i}$$
 and thus we can divide the path equally into $d_i$ subpaths
$$\Psi_{k,p,a}(n);\;\;\;a=1,...,d_i$$ each representing the word
  $$x_{b_{k-1}+p}^{n}.$$
Now we pick a vertex $u_{k,p,a}$ in $\Psi_{k,p,a}(n)$ for each $k=1,2, p=1,...., b_k, a=1,...,d_i$ as
follows
\newline
-- if $(k,p)\ne (2, b_2)$, then $u_{k,p,a}$ is the middle vertex of $\Psi_{k,p,a}(n)$ (recall that $n$
is even);
\newline
-- if $(k,p)=(2,b_2)$, then $u_{2, b_2,a}$ is a vertex around the middle vertex of
 $\Psi_{k,p,a}(n)$ which is the initial vertex of a $x_1$-edge.

By Note \ref{large power},
  for each $k=1,2$ and $p=1,...,b_k$, the set
  of $d_i$ points
   $$\{u_{k,p,a}; a=1,...,d_i\}$$ is contained in the set of $d_{i,k}$
   paths
 $$\{\Phi_{k,p,q}(n); q=1,..., d_{i,k}\}.$$

Now we  cut ${\cal{G}}_{3}(n)^f$ at the vertices   $u_{k,p,a}$,
 $k=1,2, p=1,..., b_k, a=1,..., d_i$,
that is, we
 form  a cut graph  ${\cal{G}}_{3}(n)^f_c={\cal{G}}_{3}(n)^f
 \setminus \{u_{k,p,a},
 k=1,2, p=1,..., b_k, a=1,..., d_i\}$, whose
 vertex set is obtained from
 the vertex set of ${\cal{G}}_{3}(n)^f$ by replacing each
 $u_{k,p,a}$ with a pair
 of vertices  $u_{k,p,a}^{\pm}$ (where $u_{k,p,a}^{+}$ is the  terminal vertex and
 $u_{k,p,a}^{-}$ is the initial  vertex).

 Now we take $d_i$ copies of the graph $\Omega$
   shown in Figure \ref{fg1} (a), which we denote by
   $\Omega_a$, $a=1,...,d_i$.
  For each fixed $a=1,...,d_i$, we identify the vertex set
  $$\{u_{k,p,a}^{\pm}, k=1,2, p=1,...,b_k\}$$ of ${\cal{G}}_{3}(n)^f_c$
 with  the  vertices  of  $\Omega_a$ as follows:\\
-- if $(k,p)\ne (2, b_2)$, identify $u_{k,p,a}^+$ with the left vertex of $\O_a$ if $b_{k-1}+p$ is odd
and to the right vertex if $b_{k-1}+p$ is even, and identify $u_{k,p,a}^-$ with the right vertex of
$\O_a$ if $b_{k-1}+p$ is odd and to the left vertex if $b_{k-1}+p$ is even,\\ -- identify
$u_{2,b_2,a}^+$ with the left vertex of $\O_a$  and identify $u_{2,b_2,a}^-$ with the right vertex of
$\O_a$.

\begin{figure}[!ht]
{\epsfxsize=4.7in \centerline{\epsfbox{fg2.ai}}\hspace{10mm}} \caption{}\label{fg2}
\end{figure}

 The resulting graph is not
folded, but becomes folded graph after the following obvious folding operation around each inserted
$\O_a$: \\ -- fold the subpath $x_{n_{i,j}-1}a_{1}b_{1}a_{1}^{-1}b_{1}^{-1}\cdots a_{g}
b_{g}a_{g}^{-1}b_{g}^{-1}$  whose terminal vertex is the vertex $u_{2,b_2,a}^+$ with the loops of $\O_a$
at the left vertex of $\O_a$ and then with the $x_{n_{i,j}-1}$-edge of ${\cal G}_{4}(n)$ whose terminal
vertex is the left vertex of $\O_a$, and\\ -- fold the two $x_1$-edges whose initial vertices are the
right vertex of $\O_a$.\\ The resulting folded graph, denoted  ${\cal G}_{4}(n)$,  around the inserted
$\O_a$ is shown in Figure \ref{fg2}.
 By  our construction we see that ${\cal G}_{4}(n)$ is a folded,
 $L$-labeled, directed graph, with  no
 $x_{b}$-loops, with each of the words
 $w_{1}^*,...,w_{\ell}^*$ still representable
 by a loop based at $v_{0}$, and
 with each of the words $y_{1},..., y_{r}$
 still representable by a non-closed path based at $v_{0}$.
 Also we see that the graph ${\cal G}_{4}(n)$ contains
 loops based $v_{0}$
 representing the words $$z_{k,p,q}(n+1)^*, \;\;
 k=1,2, p=1,...,b_k, q=1,...,d_{i,k}-1.$$

The graph  ${\cal G}_{4}(n)$ is not $L$-regular yet since it does not contain any $x_{b}$-loops. So it
must contain a missing label. Let $x\in L$ be a missing label at a vertex $v$ of ${\cal G}_{4}(n)$. Let
$\a$ be a finite directed graph consisting of a single path of edges all labeled with $x$, as shown in
Figure \ref{fg3}. We identify the left end vertex of $\a$ to the vertex $v$ of ${\cal G}_{4}(n)$. The
resulting
 graph ${\cal G}_{5}(n)$ is obviously still folded, contains
 ${\cal G}_{4}(n)$ as an embedded subgraph, and contains
 no $x_{b}$-loops for any $b=1,...,n_{i,j}$.
 By choosing a long enough path $\alpha$,
 we may assume that the number of vertices of ${\cal G}_{5}(n)$
 is bigger than $(n+1) d_i+1$.

\begin{figure}[!ht]
{\epsfxsize=4in \centerline{\epsfbox{fg3.ai}}\hspace{10mm}} \caption{}\label{fg3}
\end{figure}

Now by Theorem  \ref{regular},
 we can obtain an $L$-regular graph ${\cal G}_{6}(n)$ such that
 \newline
 (1)  ${\cal G}_{5}(n)$ is an embedded subgraph of ${\cal G}_{6}(n)$;
 thus in particular  in ${\cal G}_{6}(n)$
 each of the words
 $w_{1}^*,...,w_{\ell}^*$, $z_{k.p.q}(n+1)^*$,
$k=1,2, p=1,...,b_k$, $q=1,...,d_{i,k}-1$ is representable
 by a loop based at $v_{0}$, and each of the words $y_{1},..., y_{r}$
 is representable by a non-closed path based at $v_{0}$;
 \newline
 (2) ${\cal G}_{6}(n)$ contains no
 loops representing the word $x_{b}^d$ for any $b=1,...,n_{i,j}$,
  $d=1,...,m_*-1$, where $m_*$ is the number of
 vertices of ${\cal G}_{6}(n)$ (note that $m^*$ depends on $i$ and $j$).

  Note that $m_*$ is some integer larger than $(n+1)d_i+1$.
 Let $N_{i,j}=m_*-(d_i-1)(n+1)-1$. Then $N_{i,j}>(n+1)$.

 During the transformation from ${\cal G}_{4}(n)$ to ${\cal G}_{6}(n)$,
  the subgraph of ${\cal G}_{4}(n)$  consisting of the
  edges which intersect the subgraph
 $\Omega_a$ (for each fixed $a=1, ..., d_i$) remained unchanged
 since ${\cal G}_{4}(n)$ was locally $L$-regular
 already at the two vertices of
 $\O_a$.
 Now we replace $\O_a$, for each of $$a=1, ..., d_i-1,$$ by
 a graph $\O_{a}(N_{i,j}-n+1)$ which is similar to $\Omega_a$ but with $N_{i,j}-n+1 \geq 3$
 vertices (Figure \ref{fg4}
 illustrates such a graph with four vertices).
 Then the resulting graph, which we denote by  ${\cal G}(N_{i,j})$, has the following  properties.
\newline
 (1)  ${\cal G}(N_{i,j})$ is  $L$-regular;\newline
 (2) each of the words $y_{1},...,y_{r}$ is still representable by a
non-closed path based at $v_{0}$ in  ${\cal G}(N_{i,j})$,\newline
 (3) each of the words $w_{1}^*,..., w_{\ell}^*$ is still representable
 by a loop based at $v_{0}$ in  ${\cal G}(N_{i,j})$, \newline
 (4)   ${\cal G}(N_{i,j})$ contains no
 loops representing the word $x_{b}^d$ for each $b=1,...,n_{i,j}$ and
 each $d=1,...,m_\#-1$,
 where $m_\#$ is the number of vertices of  ${\cal G}(N_{i,j})$,\newline
 (5)  ${\cal G}(N_{i,j})$  contains a closed loop based at $v_{0}$
representing the word $z_{k,p,q}(N_{i,j})^*$, for each $k=1,2, p=1,...,b_k, q=1,...,d_i-1$,
and \newline
 (6) $m_\#$, the number of vertices of
 ${\cal G}(N_{i,j})$, is equal to $ N_{i,j} d_i+1$. \newline
 Properties (1)-(5) are obvious by the construction, while property (6)
follows by a simple calculation. Indeed
\begin{eqnarray*}
m_\# &=&  m_* + (N_{i,j}-n+1-2)(d_i-1)\\
 &=&[N_{i,j}  + (d_i-1)(n+1)+1] + (N_{i,j} -(n+1))(d_i-1)\\
 &=& N_{i,j} d_i+1.
\end{eqnarray*}

\begin{figure}[!ht]
{\epsfxsize=3in \centerline{\epsfbox{fg4.ai}}\hspace{10mm}} \caption{}\label{fg4}
\end{figure}

Now for each integer $N_*\geq N_{i,j}$ we construct
 a finite, connected, $L$-labeled, directed graph
${\cal G}(N_*)$
 (with a fixed base vertex
$v_0$) with the properties (0)-(5) listed in Theorem \ref{large}. In the graph ${\cal G}(N_{i,j})$
above, for each $a=1,...,d_i-1$, replace the   subgraph $\O_a(N_{i,j}-n+1)$ by the graph
$\O_a(N_*-n+1)$,  and replace subgraph $\O_{d_i}$ by the graph $\O_{d_i}(N_*-N_{i,j}+2)$. The resulting
graph is ${\cal G}(N_*)$.

\noindent {\bf Subcase (b)} $n_{i,j}>1$ is odd.

We modify the graph ${\cal G}_{3}(n)^f$
 as follows.
Besides the vertices $u_{k,p,a}$ we have chosen before, we choose, for each $a=1,...,d_i$, a vertex
$u_{2,b_2,a}^{\prime}$ in the directed subpath $\Psi_{2,b_2,a}$ such that \\-- $u_{2,b_2,a}^{\prime}$ is
 the initial vertex of an edge with label $x_2$,\\
 -- $u_{2,b_2,a}'$ appears after the vertex $u_{2,b_2,a}$
 in the directed subpath $\Psi_{2,b_2,a}$,\\
 -- there are precisely five
 edges with label $x_1$
between $u_{2,b_2,a}$ and $u_{2, b_2, a}^{\prime}$ in the directed subpath $\Psi_{2,b_2,a}$ (this
is possible as $n$ is large).

Again as $n$ is large, the set of $d_i$ vertices
   $\{u_{2,b_2,a}'; a=1,...,d_i\}$ is contained in the set of $d_{i,k}$
   paths
$\{\Phi_{2,b_2,q}(n); q=1,..., d_{i,k}\}$ (cf. Note \ref{large power}).

 Now cut ${\cal{G}}_{3}(n)^f$ at the vertices
 $\{u_{k,p,a}, k=1,2, p=1,...,b_k, a=1,...,d_i\}$, and  $\{u_{2,b_2,s}', a=1,...,d_i\}$,
 and for each $a=1,...,d_i$, insert the graph $\O_a$, which is a copy of the graph $\Omega$
   shown in Figure \ref{fg1} (b).
 That is, we\\
 (1) Form a cut graph  ${\cal{G}}_{3}(n)^f_c={\cal{G}}_{3}(n)^f\setminus
 \{u_{k,p,a},u_{2,b_2,a}',  k=1,2, p=1,...,b_k, a=1,...,d_i\}$,
 defined  as in Subcase (a), with obvious modifications,
 i.e. we have similarly defined pairs of vertices
 $u_{k,p,a}^{\pm}$, $u_{2,b_2,a}^{'\pm}$ for
 ${\cal{G}}_{3}(n)^f_c$ such that
 if each such $\pm$ pair of vertices  are identified,
 then the resulting graph is the original ${\cal{G}}_{3}(n)^f$.
\\(2) For each fixed $a=1,...,d_i$, we identify the
 vertex set $\{u_{k,p,a}^{\pm}, u_{2,b_2,a}^{'\pm}, k=1,2, p=1,...,b_k\}$ of
 ${\cal{G}}_{3}(n)^f_c$
 with  the left-most and right-most vertices  of  $\Omega_a$ as follows:\\
-- if $(k,p)\ne (2,b_2)$, and $(k, p)=(1,1)$ or
$b_{k-1}+p$ is even, then
 identify $u_{k,p,a}^+$ with the left-most vertex of
$\O_a$ and $u_{k,p,a}^-$ with the right-most vertex.\\
-- if $(k,p)\ne (2,b_2)$, $(k,p)\ne (1,1)$ and $b_{k-1}+p$ is odd, then identify $u_{k,p,a}^+$ with the
right-most vertex of $\O_a$ and $u_{k,p,a}^-$
 with the left-most vertex,\\
-- identify $u_{2,b_2, a}^+$ with the left-most vertex of $\O_a$ and identify $u_{2,b_2,a}^-$ with the
right-most vertex of $\O_a$,
\\
-- identify $u_{2,b_2, a}^{'+}$ with the left-most vertex of $\O_a$ and identify $u_{2,b_2,a}^{'-}$ with
the right-most vertex of $\O_a$.

 The resulting graph is not
folded, but  becomes folded graph after the following folding operations are performed around each
inserted  $\O_a$: \\ -- fold the path $x_{n_{i,j}-1}a_{1}b_{1}a_{1}^{-1}b_{1}^{-1}\cdots a_{g}
b_{g}a_{g}^{-1}b_{g}^{-1}$  whose terminal vertex is  the vertex $u_{2,b_2,a}^+$ with the loops of
$\O_a$ at the left-most vertex of $\O_a$ and then with the $x_{n_{i,j}-1}$-edge of ${\cal G}_{4}(n)$
whose terminal  vertex is the left-most vertex of $\O_a$,
\\ -- fold the two $x_{1}$-edges whose initial vertices are the right-most
 vertex of $\O_a$,
\\ -- fold the two $x_{1}$-edges whose terminal vertices are the left-most
vertex of $\O_a$,
\\ -- fold the two $x_{2}$-edges whose initial vertices are the right-most
vertex of $\O_a$.\\ The resulting folded graph ${\cal G}_{4}(n)^f$ around the inserted $\O_a$ is shown
in Figure \ref{fg5}.
 By  our construction we see that ${\cal G}_{4}(n)^f$ is a folded,
 $L$-labeled, directed graph, with  no
 $x_{b}$-loops, with each of the words
 $w_{1}^*,...,w_{\ell}^*$ still representable
 by a loop based at $v_{0}$, and
 with each of the words $y_{1},..., y_{r}$
 still representable by a non-closed path based at $v_{0}$.
 Also we see that the graph ${\cal G}_{4}(n)$ contains
 loops based $v_{0}$
 representing the words $z_{k,p,q}(n+2)^*$, for
 all $k=1,2, p=1,...,b_k, q=1,...,d_{i,k}-1$.

\begin{figure}[!ht]
{\epsfxsize=5.5in \centerline{\epsfbox{fg5.ai}}\hspace{10mm}} \caption{}\label{fg5}
\end{figure}

We then define ${\cal G}_{5}(n)$ and  ${\cal G}_{6}(n)$ in a similar manner as in
Subcase (a); here we
may assume that
 ${\cal G}_{5}(n)$ has at least $(d_i-1)(n+2)-1$ vertices.
 Let $m_*$ be the number of vertices of
${\cal G}_{6}(n)$, and let $N_{i,j}= m_*-(d_i-1)(n+2)-1$.
  To form  ${\cal G}(N_{i,j})$, we  replace each
 subgraph $\Omega_a$, $a=1,...,d_i-1$  in ${\cal G}_{6}(n)$
 with a graph $\O_a(1+N_{i,j}- n)$ similar to Figure \ref{fg1}(b) but with
 $1+N_{i,j}- n$ vertices.
 In the current case, we need $1+N_{i,j}-n$ to be an odd integer
 in order for the construction to work.
 This is made possible by the following

\begin{lem}
$N_{i,j} - n$ is even, i.e. $N_{i,j}$ is even (since
we have chosen $n$ to be even (see Note \ref{n even})).
\end{lem}

The proof of this lemma is similar to that of \cite[Lemma 11.3]{MZ},
 noticing in the current case the
number $d_i$ is even for each $i=1,2$ by Notation \ref{dik}.

 The rest of the argument proceeds by obvious analogy with
 the Subcase (a).
 That is, the graph ${\cal G}(N_{i,j})$ is  a graph with the
  properties listed as (1)-(6) in Subcase (a).
 Indeed, Properties (1)-(5) are immediate.
 To verify Property (6), we let $m_\#$ be the number of
 vertices of ${\cal G}(N_{i,j})$, and then we have:
\begin{eqnarray*}
m_\#&=& m_*+ (1+N_{i,j}-n-3)(d_i-1)\\
 &=&  N_{i,j} + (d_i-1)(n+2) +1+ (N_{i,j}-n-2)(d_i-1)\\
 &=&  N_{i,j}d_i+1.
\end{eqnarray*}

Now for each even integer $N_*\geq N_{i,j}$, the graph ${\cal G}(N_*)$ required by Theorem \ref{large}
is obtained from  the graph ${\cal G}(N_{i,j})$ by replacing each subgraph $\O_a(N_{i,j}-n+1)$,
$a=1,...,d_i-1$,  by the graph $\O_a(N_*-n+1)$, and replacing the subgraph $\O_{d_i}$ by the graph
$\O_{d_i}(N_*-N_{i,j}+3)$.

\noindent{\bf Proof of Theorem \ref{each large n} in Case 2}.

The proof is similar to that of Case 1 and much simpler notationally, and we shall be very brief.
 In
this case $n_{i,j}=1$, i.e. the surface $S_{i,j}$ has only one boundary component, which we denote by
$\b$ and  may assume lying in $T_1$.  $\b$ has $d_{i,1}$ intersection points with $\p S_{i_*}$, which we
denote by $t_q$, $q=1,...,d_{i,1}$. Similarly as in Case 1, we define the points $b_q, q=1,...,d_{i,1}$
in $\p_p J_{i,j}$. The group $\pi_1(S_{i,j}^-,t_1)$ has a set of generators
$$L=\{a_1, b_1,...,a_g, b_g\}$$
($g$ must be larger than $0$) such that $$x_1=a_1 b_1 a_1^{-1}b_1^{-1}\cdots a_g b_g a_g^{-1}b_g^{-1} $$
is an embedded loop which is homotopic to $\b$. As in Case 1, we similarly  define the elements
$w_1,...,w_{\ell}$, the element $y_1,...,y_r$, and the elements $z_q(n)$, $q=1,...,d_{i,1}-1$, and we
reduce the proof of Theorem \ref{each large n} in Case 2 to the proof of the following theorem which is
an analogue of Theorem \ref{large}.

\begin{thm}\label{large2} There is a positive even integer $N_{i,j}$ such that for
each even integer $N_*\geq N_{i,j}$ there is a finite, connected, $L$-labeled, directed graph ${\cal
G}(N_*)$
 (with a fixed base vertex
$v_0$) with the following properties:
\newline
  (0)  ${\cal G}(N_*)$ is  $L$-regular;\newline
(1) The number of vertices of  ${\cal G}(N_*)$ is
 $m_{i}=N_*d_i+1$; \newline
 (2) Each of the words $w_1^*,..., w_\ell^*$ is  representable
 by a loop, based at $v_0$, in ${\cal G}(N_*)$; \newline
(3) ${\cal G}(N_*)$  contains a closed loop, based at $v_0$, representing the word $z_{q}(N_*)^*$, for
each $q=1,...,d_{i,1}-1$;\newline
 (4)  ${\cal G}(N_*)$ contains no
 loop representing the word $x_1^d$ for any
 $d=1,...,m_{i}-1$;\newline
  (5) each of the words $y_1,...,y_r$ is  representable by a
non-closed path, based at  $v_0$, in ${\cal G}(N_*)$.
\end{thm}

To prove this theorem, we construct, similar as in Case 1,
 the analogue graph ${\cal G}_3(n)^f$ and its subgraphs
 $\Phi_q(n), q=1,...,d_{i,1}$, $\Psi_a(n), a=1,...,d_i$, with similar properties.
 We modify the graph ${\cal G}_{3}(n)^f$
 as follows.
For each of $a=1,...,d_i$, we pick a pair  vertices $\{u_a,
u_{a}'\}$ in the path $\Psi_a(n)$ as follows:\\
--$u_a$ is closed to the middle vertex of $\Psi_a(n)$;\\
  --$u_a$  is the terminal vertex of an edge
with label $a_{1}$; and \\
 --$u_{a}^{\prime}$  is the terminal vertex of an edge
with label $b_{1}$  which appears after the vertex $u_{a}$;\\
--there are precisely five  edges with label $b_{1}$ between $u_{a}$ and $u_{a}^{\prime}$ in the path
$\Psi_{a}(n)$.

We may assume that the set  $$\{u_a, u_a'; a=1,...,d_i\}$$  is contained in the set
 $$\{\Phi_q(n); q=1,...,d_{i,1}\}.$$

\begin{figure}[!ht]
{\epsfxsize=3in \centerline{\epsfbox{fg6.ai}}\hspace{10mm}} \caption{}\label{fg6}
\end{figure}

 Now cut the graph ${\cal{G}}_{3}(n)^f$ at all the pairs
 of vertices $\{u_{a}, u_{a}'\}$, $a=1,...,d_i$,
and for each $a$, insert the graph $\O_a$-- which is a copy of the graph $\Omega$
   shown in Figure \ref{fg6} -- as follows.
 Form a cut graph  ${\cal{G}}_{3}(n)^f_c={\cal{G}}_{a}(n)^f\setminus
 \{u_{a}, u_{a}'; a=1,...,d_i\}$,
 and let
 $u_{a}^{\pm}$, $u_{a}^{'\pm}$ be the corresponding vertices
 for
 ${\cal{G}}_{3}(n)^f_c$.
 For each fixed $a=1,...,d_i$, we identify the
 vertex $u_{a}^+$ with the left-most vertex of $\O_a$, identify
$u_{a}^-$ with the right-most vertex of $\O_a$,
 identify $u_{a}^{'+}$ with the right-most vertex of $\O_a$
and identify $u_{a}^{'-}$ with the left-most  vertex of $\O_a$.

 The resulting graph is not
folded, but  becomes folded graph after a single folding operation around each inserted  $\O_a$:
 fold the two
$a_{1}$-edges whose terminal  vertices are the right-most vertex of $\O_a$. The resulting folded graph
${\cal G}_{4}(n)^f$ around the inserted $\O_a$ is shown in Figure \ref{fg7}.
 By  our construction we see that ${\cal G}_{4}(n)^f$ is a folded
 $L$-labeled directed graph, with  no
 $x_1$-loops, with each of the words
 $w_{1}^*,...,w_{\ell}^*$ still representable
 by a loop based at $v_{0}$, and
 with each of the words $y_{1},..., y_{r}$
 still representable by a non-closed path based at $v_{0}$.
 Also we see that the graph ${\cal G}_{4}(n)$ contains
 loops based at $v_{0}$
 representing the words $z_{q}(n+4)^*$, for
 all $q=1,...,d_{i,1}-1$.

\begin{figure}[!ht]
{\epsfxsize=5.5in \centerline{\epsfbox{fg7.ai}}\hspace{10mm}} \caption{}\label{fg7}
\end{figure}

As in Case 1, we get ${\cal G}_{5}(n)$ and ${\cal G}_{6}(n)$. In the current case, $N_{i,j}= m_*-
(d_i-1)(n+4)-1$, which is   larger than $n+4$, where $m_*$ is the number of vertices of ${\cal G}_6(n)$.
  To form  ${\cal G}(N_{i,j})$, we replace the left half
 (with three vertices) of
   $\Omega_a$, for each  $a=1,...,d_i-1$,
 with a graph $\O_a(N_{i,j}- n-1)$ which is similar to Figure \ref{fg6} but with
 $N_{i,j}- n-1$ vertices.
 In the current case, we also need $N_{i,j}$ to be an even integer
 in order for the construction to work.
 This is true, and  can be proved as in Subcase (b) of Case 1.
It is easy to see that ${\cal G}(N_{i,j})$ has all the Properties (0)-(5) listed in Theorem \ref{large2}
(when $N_*=N_{i,j}$). For instance to check Property (1), we have:
\begin{eqnarray*}
m_i &=& m_*+(N_{i,j}-n-1-3)(d_i-1)\\
    &=& N_{i,j}+(d_i-1)(n+4)+1+ (N_{i,j}-n-4)(d_i-1)\\
 &=& N_{i,j}d_i+1
\end{eqnarray*}
To show that Theorem \ref{large2} holds for any even integer $N_*\geq N_{i,j}$, we simply let ${\cal
G}(N_*)$ be the graph  obtained from  the graph ${\cal G}(N_{i,j})$ by replacing each subgraph
$\O_a(N_{i,j}-n-1)$, $a=1,...,d_i-1$,  by the graph $\O_a(N_*-n-1)$, and replacing the subgraph
$\O_{d_i}$ by the graph $\O_{d_i}(N_*-N_{i,j}+5)$.

\section{The final assembly}\label{hsmf}

Fix an even integer $N_*$ satisfying Corollary \ref{indep of j}
(later on we may need $N_*$ to have been chosen large enough).
Then as given in Corollary \ref{indep of j}, for each $i=1,2$ and $j=1,.., n_i$, the
manifold $Y_{i,j}^-=S_{i,j}^-\times I$
has an $m_i=N_* d_i+1$ fold cover
$\breve Y_{i,j}^-=\breve S_{i,j}^-\times I$ such that
$|\p_p   \breve Y_{i,j}^-|= |\p_p  Y_{i,j}^-|$, i.e.
 each component of $\p_p  \breve Y_{i,j}^-$ is
 an $m_i$ fold cyclic cover of a component of
 $\p_p   Y_{i,j}^-$.
 Moreover the map $g_{i,j}: (J_{i,j}^-, \p_p J_{i,j}^-)
 \ra (Y_{i,j}^-, \p_p Y_{i,j}^-)$ lifts to an embedding
 $\breve g_{i,j}: (J_{i,j}^-, \p_p J_{i,j}^-)
 \ra (\breve Y_{i,j}^-, \p_p \breve Y_{i,j}^-)$ such that
if $\breve A$ is a component of $\p_p\breve Y_{i,j}^-$ then
the components of  $\breve g_{i,j}(\p_p J_{i,j}^-)\cap \breve A$
are evenly distributed along $\breve A$.
More precisely if
$S_{i,j}^-$, for instance,
 is the surface given in the proof of Theorem \ref{each large n} in Case 1,
 then with the notations given there, we may index the boundary components of $\breve S_{i,j}^-$
 as $\breve\b_{k,p}$, $k=1,2,p=1,...,b_k$,
 so that each $\breve\b_{k,p}$ is an $m_i$ fold cyclic cover of
 $\b_{k,p}$ and the points
 $\{\breve g_{i,j}(b_{k,p,q}); q=1,...., d_{i,k}\}$
 divide $\breve \b_{k,p}$ into $d_{i,k}$ segments each
 having wrapping number
 $N_* d_i/d_{i,k}$. Also note that the $d_{i,k}$ points
 $\{\breve g_{i,j}(b_{k,p,q}); q=1,...., d_{i,k}\}$ are mapped to the $d_{i,k}$ points
 $\{t_{k,p,q}; q=1,...., d_{i,k}\}$ respectively under the covering
  map $\breve\b_{k,p}\ra \b_{k,p}$.
As $N_*$ can be assumed to be arbitrarily large, we may assume that
 the wrapping number $N_* d_i/d_{i,k}$ be as large as needed
 for each $i=1,2$ and $k=1,...,m$.

 Also recall
that  $(K_{i,j}^-, \p_p K_{i,j}^-)$ is properly embedded in the pair $(J_{i,j}^-, \p_p J_{i,j}^-)$,
with a relative $R$-collared neighborhood. It follows that the pair $(\breve g_{i,j}(K_{i,j}^-), \breve
g_{i,j}(\p_p K_{i,j}^-))$ has a relative $R$-collared neighborhood in $(\breve Y_{i,j}^-, \p_p\breve
Y_{i,j}^-)$. Also  $K_1^-=\cup_{j=1}^{n_1} K_{1,j}$ and $K_2^-=\cup_{j=1}^{n_2} K_{2,j}$ are isometric
under the isometry $h:K_1\ra K_2$. Now let $\breve Y^-$ be the union of $\breve Y_1^-=\cup_{j=1}^{n_1}
Y_{1,j}^-$ and $\breve Y_2^-=\cup_{j=1}^{n_2} Y_{2,j}^-$ with
 $\cup_{j=1}^{n_1} (\breve g_{1,j}(K_{1,j}^-), \breve
g_{1,j}(\p_p K_{1,j}^-))$ and  $\cup_{j=1}^{n_2} (\breve g_{2,j}(K_{2,j}^-), \breve g_{2,j}(\p_p
K_{2,j}^-))$ identified by the corresponding isometry and let $(U^-, \p_p U^-)$ be the identification
space of  $\cup_{j=1}^{n_1} (\breve g_{1,j}(K_{1,j}^-), \breve g_{1,j}(\p_p K_{1,j}^-))$ and
$\cup_{j=1}^{n_2} (\breve g_{2,j}(K_{2,j}^-), \breve g_{2,j}(\p_p K_{2,j}^-))$
 in $\breve Y^-$. Then $\breve Y^-$ is a connected metric space,
 with a path metric induced from the metrics on $\breve Y_1^-$ and $\breve Y_2^-$.
 There is an induced local isometry $f:\breve{Y}^- \ra M$ which extends
 the local isometry
$\breve{Y}_{i,j}^- \ra Y_{i,j}^-\ra M$ for each $i,j$.

 Define the parabolic boundary, $\p_p\breve Y^-$, of $\breve
Y^-$ to be the union of $\p_p\breve Y_1^-$ and $\p_p\breve Y_2^-$, with
 $\cup_{j=1}^{n_1}  \breve
g_{1,j}(\p_p K_{1,j}^-)$ and  $\cup_{j=1}^{n_2}\breve g_{2,j}(\p_p K_{2,j}^-)$  identified by the
isometry. Then $(U^-, \p_p U^-)$ has a relative $R$-collared neighborhood in $(\breve Y^-, \p_p\breve
Y^-)$.

Now for each $k=1,...,m$, let $\breve C_K$ be the cover of the $k$-th cusp $C_k$ of $M$ corresponding to
the subgroup of $\pi_1(C_k)$ generated by the $m_1$-th power of a component of $\p_k S_1^-$ and the
$m_2$-th power of a component of $\p_k S_2^-$. Then  each oriented component, say $\b$,  of $\p_k S_i^-$
has its inverse image in $\p\breve C_k$, denoted $\breve \b$, a connected oriented circle. So we may and
shall identify $\breve\b$ with the oriented component of $\p \breve S_i^-$ which covers $\b$. This way
we embed naturally all components of $\p\breve S_i^-$ into $\p \breve C=\cup _{k=1}^m \p \breve C_k$,
for each $i=1,2$. We denote by $\p_k \breve S_i^-$ those components of $\p \breve S_i^-$ which are
embedded in $\p \breve C_k$. Then we have  $|\p_k \breve S_i^-|=|\p_k S_i|$, and the components of $\p_k
\breve S_i^-$  are far apart from each other in $\p \breve C_k$, for each  $i=1,2$ and $k=1,...,m$. So
we may and shall embed the corresponding components of $\p_p\breve Y_i^-$ into $\p \breve C$ along $\p
\breve S_i^-$, for each $i=1,2$. After such identification, we get a connected hyperbolic $3$-manifold
$$\breve Y^-\cup(\cup_{k=1}^m \breve C_k)$$ with $m$ rank two cusps and with a local isometry into $M$.

As in \cite[Section 13]{MZ} we construct the thickening $\bar U^-$ of $U^-$ so that $\p_p \bar U$ is embedded in $\p \breve C$
(Note that  each component of $\bar U^-$ is a
handlebody, with a similar proof as that of  \cite[Lemma  13.2]{MZ})
and let
$$Y^-=\breve Y_1^-\cup \bar U^-\cup \breve Y_2^-.$$
Then $Y^-$  is a connected, compact, hyperbolic $3$-manifold, locally convex everywhere except on its
parabolic boundary $\p_p Y^-=\p_p \breve Y_1^-\cup \p_p \bar U^-\cup \p_p\breve Y_2^-$.
 The complement of $\p_p(Y^-)$ in $\p \breve C$ is
  a set of ``round-cornered parallelograms'' with very long sides
 in $\p \breve
C$. As in \cite[Section 13]{MZ}, we scoop out from  $\breve C=\cup_{k=1}^m C_k$ the convex half balls
based on these round-cornered  Euclidean parallelograms and denote the resulting cusps by $\cup_{k=1}^m
\breve C_k^0$. Then
$$Y=Y^-\cup (\cup_{k=1}^m\breve C_k^0)$$
is a connected, metrically complete, convex hyperbolic $3$-manifold, with a local isometry $f$
into $M$.
Thus the local isometry $f$ induces an injection of $\pi_1(Y)$ into $\pi_1(M)$
by \cite[Lemma 4.2]{MZ}.

Each boundary component of $Y$ provides a
Quasi-Fuchsian surface in $M$.
To prove this claim, it suffices to show,
with a similar reason as that given in \cite[Section 13]{MZ},
that   every Dehn filling of $Y$ along its cusps
gives a $3$-manifold with incompressible boundary.

Let $Y(\a_1,...,\a_m)$ be any  Dehn filling of $Y$ with slopes $\a_1,...,\a_m$.
 Then  $Y(\a_1,...,\a_m)$
is an $HS$-manifold (see  \cite[Section 12]{MZ} for its definition).
 The handlebody part
$H$  of $Y(\a_1,...,\a_m)$ is
 $\bar U^-\cup (\cup_{k=1}^m \breve
C_k^0(\a_k))$ (which may have several components in the current case but each has genus at least
two) and the $S\times
I$ part of $Y(\a_1,...,\a_m)$ is $Y(\a_1,...,\a_m)\setminus H=Y^-\setminus \bar U^-$.
This $HS$ decomposition of $Y(\a_1,...,\a_m)$ satisfies the conditions of \cite[Lemma 12.1]{MZ}
and thus $Y(\a_1,...,\a_m)$ has incompressible boundary by that lemma.
The proof of this last claim is similar to that of
\cite[Lemma 13.6]{MZ}, for which we only need to note the following:
\newline
(1)
With a similar proof as that of  \cite[Lemma  13.5]{MZ} we have that each component of
$Y^-\setminus \bar U^-$ is not simply connected.
\newline
(2) $|\p_k \breve S^-_i|=|\p_k S_i^-|\geq 2$ for each $i=1,2, k=1,...,m$, by
Condition  \ref{at least two}.

The proof of Theorem \ref{main} is now finished.

\end{document}